\documentclass[12pt]{article}
\newtheorem{theo}{{\bfseries Theorem}}[section]
\newtheorem{prop}[theo]{{\bfseries Proposition}}
\newtheorem{lem}[theo]{{\bfseries Lemma}}
\newtheorem{cor}[theo]{{\bfseries Corollary}}
\newtheorem{df}[theo]{{\bfseries Definition}}

\newtheorem{ques}[theo]{{\bfseries Question}}

\def \R {\mathbb R}

\def \I {\mathcal I}
\def \J {\mathcal J}

\def \a {\alpha }
\def \b {\beta}
\def \e {\epsilon}
\def \d {\delta}
\def \g {\gamma}

\def \ee {{\mathbf e}}
\def \pp {{\mathbf p}}
\def \qq {{\mathbf q}}
\def \vv {{\mathbf v}}
\def \MM {{\mathbf M}}
\def \tp {\tilde {\mathbf p}}
\def \tq {\tilde \tilde{\mathbf q}}
\def \SX {{\mathbf S_X}}
\def \SY {{\mathbf S_Y}}
\def \one {{\mathbf 1}}

\usepackage{amssymb,latexsym}
\usepackage[mathcal]{eucal}
\usepackage{amscd}
\usepackage{amsmath}
\usepackage{graphicx}
\usepackage{amsfonts}
\usepackage{amscd}

\numberwithin{equation}{section}

\renewcommand{\arraystretch}{1.5}

\begin{document}

\title{\bfseries  The Iterated Prisoner's Dilemma: Good Strategies and Their Dynamics}
\vspace{1cm}
\author{Ethan Akin\\
    Mathematics Department \\
    The City College \\
    137 Street and Convent Avenue \\
    New York City, NY 10031, USA }

    \vspace{.5cm}
\date{July, 2013}

 \maketitle

\begin{abstract}  For the iterated Prisoner's Dilemma, there exist  Markov strategies
which solve the problem when we restrict attention to the long term
average payoff. When used by both players these assure the
cooperative payoff for each of them. Neither player can benefit by
moving unilaterally any other strategy, i.e. these are Nash
equilibria. In addition, if a player uses instead an alternative
which decreases the opponent's payoff below the cooperative level,
then his own payoff is decreased as well. Thus, if we limit
attention to the long term payoff,  these \emph{good strategies} effectively
stabilize cooperative behavior. We characterize these good strategies and analyze their role in
evolutionary dynamics.
\end{abstract} \vspace{.5cm}


\vspace{1cm}

\section{The Iterated Prisoner's Dilemma}

The \emph{Prisoner's Dilemma} is a two person game that provides a simple model of a disturbing social phenomenon.
It is a symmetric game in which each of the two players, X and Y, has a choice between
two strategies, $c$ and $d$. Thus, there are four outcomes which we list in the order: $cc, cd, dc, dd,$ where,
for example,  $cd$ is the outcome when X plays $c$ and Y plays $d$. Each then receives a payoff. The following
 $2 \times 2$ chart describes the payoff to the X player. The transpose is the Y payoff.
\renewcommand{\arraystretch}{1.5}
\begin{equation}\label{1}
\begin{array}{|c||c|c|}\hline
X\backslash Y & \quad c \quad & \quad  d \quad \\ \hline \hline
c & \quad R \quad & \quad S \quad  \\ \hline
d & \quad T \quad  & \quad P \quad \\ \hline
\end{array}
\end{equation}

Alternatively, we can define the \emph{payoff vectors} for each player by
\begin{equation}\label{1a}
\SX \quad = \quad (R, S, T, P) \qquad \mbox{and} \qquad \SY \quad = \quad (R, T, S, P).
\end{equation}

Davis [6] and Straffin [17] provide clear introductory
discussions of the elements of game theory.

Either player can use a \emph{mixed  strategy},   randomizing by choosing $c$ with probability $p_c$ and
$d$ with the complementary probability $1-p_c$.

A probability \emph{distribution} $\vv$ on the set of outcomes is a non-negative
vector with unit sum, indexed by the four states.
That is, $v_i \geq 0$ for $i = 1,...,4$ and the dot product $ < \vv \cdot \one > \ = \ 1$.
For example, $v_2$ is the probability that X played $c$ and Y played $d$.
In particular, $v_1 + v_2$ is the probability X played $c$.
With respect to $\vv$ the expected payoffs to X and Y, denoted $s_X$ and $s_Y$, are
the dot products with the corresponding payoff vectors:
\begin{equation}\label{1new}
s_X \quad = \quad < \vv \cdot \SX > \qquad \mbox{and} \qquad s_Y \quad = \quad < \vv \cdot \SY >.
\end{equation}

The payoffs are assumed to satisfy
\begin{equation}\label{2}
T \ > \ R  \ > \ P  \ > \ S \qquad \mbox{and} \qquad 2R  \ > \ T + S.
\end{equation}

We will later use the following easy consequence of these inequalities.

\begin{prop}\label{newprop1} If $\vv$ is a distribution, then the associated expected payoffs to the two players,
as defined by (\ref{1new}), satisfy the following equation.
\begin{equation}\label{2new}
s_Y \ - \ s_X \quad  = \quad  (v_2 - v_3)(T - S). \hspace{2cm}
\end{equation}
So we have
$s_Y \  = \  s_X $ iff $ v_2 \  = \  v_3.$

In addition,
\begin{equation}\label{3new}
\frac{1}{2}(s_Y + s_X) \quad \leq \quad R,
\end{equation}
with equality iff $\vv \ = \ (1, 0, 0, 0)$.  Hence, the following statements are
equivalent.
\begin{itemize}
\item[(i)] $\frac{1}{2}(s_Y + s_X) \  = \  R.$
\item[(ii)] $v_1 \  = \  1$.
\item[(iii)] $s_Y \ = \ s_X \ = \ R.$
\end{itemize}
\end{prop}

{\bfseries Proof:}  Dot $\vv$ with $\SY - \SX = (0, \ T - S, \ S - T, \ 0)$ and with
$ \frac{1}{2}(\SY + \SX) = (R, \ \frac{1}{2}(T+S), \ \frac{1}{2}(T+S), \ P)$. Observe that R is the
maximum entry of the latter. \hspace{4cm}$\Box$
 \vspace{.5cm}

In the Prisoner's Dilemma, the strategy $c$ is \emph{cooperation}.  When both players cooperate
they each receive the reward for cooperation (= $R$).
The strategy $d$ is \emph{defection}.  When both players defect they each receive the punishment for defection (= $P$).
However, if one player cooperates and the other does not, then the defector receives the large temptation payoff (= $T$), while
the hapless cooperator receives the very small sucker's payoff (= $S$). The condition $2R > T + S$ says that the reward
for cooperation is larger than the players would receive by dividing equally the total payoff of a $cd$ or $dc$ outcome.
Thus, the maximum total payoff occurs uniquely at $cc$ and that location is a \emph{strict Pareto optimum}, which means that
at every other outcome at least one player does worse. The cooperative
outcome $cc$ is clearly where the players ``should" end up. If they
could negotiate a binding agreement in advance of play, they would agree to play $c$ and each receive R. However, the structure
of the game is such that, at the time of play, each chooses a strategy in ignorance of the other's choice.

This is where it gets ugly.  In game theory lingo, the strategy $d$ \emph{strictly dominates} strategy $c$.  This means that,
whatever Y's choice is, X receives a larger payoff by playing $d$ than by using $c$. In the array (\ref{1}) each number in the
$d$ row is larger than the corresponding number in the $c$ row above it. Hence, X chooses $d$, and for exactly the same reason,
Y chooses $d$. So they are driven to the $dd$ outcome with payoff P for each. Having firmly agreed to cooperate, X
hopes that Y will stick to the agreement because  X can then obtain the large payoff T by defecting. Furthermore, if he
were not to play $d$, then he risks getting S when Y defects. All the more reason to defect, as X realizes
Y is thinking the same thing.

The payoffs are often stated in money amounts or in years reduced from a prison
sentence (the original ``prisoner" version), but it is
important to understand that the payoffs are really in units of \emph{utility}.
That is, the ordering in (\ref{2}) is assumed to describe the
order of desirability of the various outcomes to each player when all the consequences of each outcome are taken into account.
Thus, if X is induced to feel guilty at the $dc$ outcome, then the payoff to X of that outcome is reduced. Adjusting the payoffs
is the classic way of stabilizing cooperative behavior.  Suppose prisoner X walks out of prison, free after defecting, having
consigned Y, who played $c$, to a 20 year sentence. Colleagues of Y might well do X some serious damage.
Anticipation of such an event
considerably reduces the desirability of $dc$ for X, perhaps to well below R. If X and Y each have
threatening friends, then it is reasonable
for each to expect that a prior agreement to play $cc$ will stand and so they each receive R.
However, in terms of utility this is
no longer a Prisoner's Dilemma.  In the  book which  originated modern game theory,
Von Neumann and Morgenstern [19], the authors
developed an axiomatic theory of utility which allows us to make sense of such arithmetic
relationships as the second inequality in (\ref{2}).  We won't consider
this further, but the reader should remember that the payoffs are numerical measurements of desirability.

This two person collapse of cooperation can be regarded as a simple model
of what Garret Hardin [7] calls \emph{the tragedy of the commons}.
 This is a similar sort of collapse of mutually beneficial cooperation on a multi-person scale.

In the search for a way to avert this tragedy,
attention has focused upon \emph{repeated play}.  X and Y play repeated rounds of the same game. For each round
the players' choices are made independently, but each is aware of all of the previous outcomes. The hope is that the threat of
future retaliation will rein in the temptation to defect in the current round.

Robert Axelrod devised a tournament in which submitted computer programs played against one another.
Each program played  a fixed, but unknown, number of rounds against each of the competing
programs, and the resulting payoffs were summed. The results are described
and analyzed in his landmark book [4].
The winning program, Tit-for-Tat, submitted by
game theorist Anatol Rapaport, consists, after initial cooperation, in playing in each round the strategy used by the opponent in the
previous round. A second tournament yielded the same winner. Axelrod extracted some interesting
rules of thumb from Tit-for-Tat and applied these to some historical examples.

At around the same time, game theory was being introduced by John Maynard Smith into biology in order to study problems
in the evolution of behavior. Maynard Smith [11] and Sigmund [13] provide good surveys of the early work.
Tournament play for games, which has been widely explored since, exactly simulates the dynamics examined in
this growing field of evolutionary game theory. However, the tournament/evolutionary viewpoint changes the problem
in a subtle way. In evolutionary game theory, what
matters is how a player is doing as compared with the competing players. Consider this with just two players and
suppose they are currently considering strategies with the same payoff to each. Comparing
outcomes, Y would reject a move to a strategy
where she does better, but which allows X to do still better than she.  That this sort of \emph{altruism} is
selected against is a major problem in the theory of evolution. However, in classical game theory
the payoffs are in utilities.  Y simply desires to obtain the highest absolute payoff.  The payoffs to her opponent
are irrelevant, except as data to predict X's choice of strategy. It is the classical problem that we will mainly consider, although
we will return to evolutionary dynamics in the last section.

I am not competent to summarize the immense literature devoted to these matters. I recommend the excellent book
length treatments of Hofbauer and Sigmund [9], Nowak [12] and Sigmund [14]. The latter two discuss
 the Markov approach which we now examine.

The choice of play for the first round is the \emph{initial play}.
A \emph{strategy}  is a choice of initial play together with what we will call a \emph{ plan}:
a choice of play, after the first round, to respond to any possible past history of outcomes
in the previous rounds. A \emph{memory-one  plan}
bases its response entirely on outcome of the previous round.
The Tit-for-Tat plan  (hereafter, just $TFT$) is an example of a \emph{ memory-one plan}.

With the outcomes listed in order as $cc, cd, dc, dd$, a memory one plan for X is a vector
$\pp = (p_1, p_2, p_3, p_4) = (p_{cc},p_{cd},p_{dc},p_{dd})$ where $p_z$
is the probability of playing c when the  outcome $z$ occurred in the previous round. If Y uses the memory-one plan
$\qq = (q_1,q_2,q_3,q_4) $ then
the  response vector is $(q_{cc},q_{cd},q_{dc},q_{dd}) = (q_1, q_3, q_2, q_4)$ and
 the successive outcomes follow a Markov chain with transition matrix
given by:
\begin{equation}\label{3}
\MM \quad =\quad \begin{pmatrix}p_1q_1 & p_1(1-q_1) & (1-p_1)q_1 & (1-p_1)(1-q_1)  \\
p_2q_3 & p_2(1-q_3) & (1-p_2)q_3 & (1-p_2)(1-q_3) \\
p_3q_2 & p_3(1-q_2) & (1-p_3)q_2 & (1-p_3)(1-q_2) \\
p_4q_4 & p_4(1-q_4) & (1-p_4)q_4 & (1-p_4)(1-q_4)
\end{pmatrix}.
\end{equation}
We use the switch in numbering from the Y plan $\qq $ to the Y response vector because switching the perspective
of the players interchanges $cd$ and $dc$. This way
the ``same" plan for X and for Y is given by the same vector.  For example, \emph{TFT} for X and for Y
is given by $\pp  = \qq = (1,0,1,0)$, but the response vector for Y is $ (1,1,0,0)$.  The plan   \emph{Repeat}
is given by $\pp   = \qq = (1,1,0,0)$ with the response vector for Y  equal to $ (1,0,1,0)$. This plan just repeats
the previous play, regardless of what the opponent did.

We describe some elementary facts about finite Markov chains, see, e.g., Chapter 2 of Karlin and Taylor [10].

A Markov matrix like $\MM$ is a non-negative matrix with row sums equal to $1$.  Thus, the vector $\one$ is a right
eigenvector with eigenvalue $1$. For such a
matrix, we can represent the associated Markov chain as movement along a directed graph with vertices the states, in this case,
$cc,cd,dc,dd$, and with a directed edge from the $i^{th}$ state $z_i$ to the $j^{th}$ state $z_j$ when
 $\MM_{ij} > 0$, that is, when
we can move from $z_i$ to $z_j$ with positive
probability. In particular, there is an edge from $z_i$ to itself iff the diagonal entry $\MM_{ii}$ is positive.

  A \emph{path} in the graph is a state sequence $z^1,...,z^n$ with $n > 1$ such that there is an
edge from $z^i$ to $z^{i+1}$ for $i = 1,...,n-1$.  A set of states $I$ is called a \emph{closed set} when
no path  that begins in $I$ can
exit $I$.  For example, the entire set of states is closed and for any  $z$ the
set of states accessible via a path that begins at $z$
is a closed set. $I$ is closed iff $\MM_{ij} = 0$ whenever $z_i \in I$ and $z_j \not\in I$.
In particular, when we restrict
the chain to a closed set $I$, the associated submatrix of $\MM$ still has row sums equal to $1$. A minimal, nonempty, closed
set of states is called a \emph{terminal set}. A state is called \emph{recurrent} when it lies in some
terminal set and \emph{transient} when it does not. The following facts are easy to check.
\begin{itemize}
\item A nonempty, closed set of states $I$ is terminal
iff for all $z_i, z_j \in I$,  there exists a path from  $z_i$ to $z_j$.
\item If $I$ is a terminal set and $z_j \in I$, then there
exists $z_i \in I$ with an edge from $z_i$ to $z_j$.
\item Distinct terminal sets are disjoint.
\item Any nonempty, closed set contains at least one terminal set.
\item  From any transient state there is a path into some terminal set.
\end{itemize}

Suppose we are given an initial distribution $\vv^1$, describing the outcome of the
first round of play. The Markov process evolves
in discrete time via the equation
\begin{equation}\label{3a}
\vv^{n+1} \quad = \quad  \vv^{n} \cdot \MM,
\end{equation}
where we regard the distributions as row vectors.

In our game context, the initial distribution is given by the
initial plays, pure or mixed, of the two players. If X uses initial probability $p_c$ and Y uses $q_c$, then
\begin{equation}\label{3aa}
\vv^1 \quad = \quad ( p_c q_c,  p_c(1-q_c), (1-p_c)q_c, (1-p_c)(1-q_c)). \hspace{1cm}
\end{equation}

Thus, $v^n_i$ is
the probability that outcome $z_i$ occurs on the $n^{th}$ round of
play. A distribution $\vv$ is \emph{stationary} when it satisfies
$\vv \MM = \vv$. That is, it is a left eigenvector with
eigenvalue $1$. From  Perron-Frobenius theory (see, e.g.,
 Appendix 2 of [10]) it follows that if $I$ is a terminal
set, then there is a unique stationary distribution $\vv$ with
$v_i > 0$ iff $i \in I$. That is, the \emph{support} of $\vv$ is exactly $I$.
In particular, if the eigenspace of $\MM$
associated with the eigenvalue $1$ is one dimensional, then there is
a unique stationary distribution, and so a unique terminal set which
is the support of the stationary distribution. The converse is also
true and any stationary distribution $\vv$ is a mixture of the $\vv_J$'s where $\vv_J$
is supported on the terminal set $J$.  This follows from the fact that any stationary distribution
$\vv$ satisfies $v_i = 0$ for all transient states $z_i$ and so is
supported on the set of recurrent states. On the recurrent states the matrix $\MM$ is block diagonal. Hence, the following are
equivalent in our $4 \times 4$ case.
\begin{itemize}
\item     There is a unique terminal set of states for the process associated with $M$.
\item     There is a unique stationary distribution vector for $M$.
\item     The matrix $M' = M - I$ has rank $3$.
\end{itemize}
We will call $\MM$ \emph{convergent} when these conditions hold. For example, when all of the probabilities of $\pp$ and $\qq$
lie strictly between 0 and 1, then all the entries of $\MM$ given by (\ref{3}) are positive and
so  the entire
set of states is the unique terminal state and  the positive matrix $\MM$ is convergent.

The sequence of the Cesaro averages $\{ \frac{1}{n} \Sigma_{i = 1}^{n} \ \vv^{i} \}$
 of the outcome distributions always
converges to some stationary distribution $\vv$. That is,
\begin{equation}\label{3ab}
Lim_{n \to \infty} \ \frac{1}{n} \Sigma_{k = 1}^{n} \ \vv^{k} \quad = \quad \vv.
\end{equation}
Hence, using the payoff vectors from (\ref{1a})
the long run average payoffs for X and Y converge to $s_X$ and $s_Y$ of (\ref{1new}) with $\vv$ this limiting stationary
distribution.

When $\MM$ is  convergent, the limit $\vv$ is the unique stationary distribution and so
the average payoffs are independent of the initial distribution.
In the non-convergent case, the long term payoffs depend on the
initial distribution. Suppose there are exactly two terminal sets,
$I$ and $J$, with stationary distribution vectors $\vv_I$ and $\vv_J$,
supported on $I$ and $J$, respectively. For any initial distribution
$\vv^1$, there are probabilities $p_I$ and $p_J = 1 - p_I$ of
entering into, and so terminating in, $I$ or $J$, respectively. In that
case, the limit of the Cesaro averages sequence for  $\{ \vv^n \}$  is given by
\begin{equation}\label{5}
\vv \quad = \qquad p_I \vv_I \ + \ p_J \vv_J, \hspace{1cm}
\end{equation}
and the limits of the average payoffs  are  given by (\ref{1new}) with this distribution $\vv$.
This extends in the obvious way when there are more terminal sets.

When Y responds to the memory-one plan $\pp$ with a memory-one plan $\qq$, we have the \emph{Markov case} as above.
We will also want to see how a memory-one plan $\pp$ for X fares against a not necessarily memory-one response by Y.
We will call such a response pattern a  \emph{general plan} to emphasize that it need not be memory-one
That is, a general plan is a choice of response, pure or mixed, for any sequence of previous outcomes.
Hereafter, unless we use the expression ``general plan'', we will assume a
 plan is  memory-one.

If Y uses a general plan, then  the sequence of Cesaro averages need not converge. We will call any  limit point of the
 sequence  \underline{an} associated \emph{limit distribution}. We will call $s_X$ and $s_Y$, given by (\ref{1new}) with
 such a limit distribution $\vv$, the \emph{expected payoffs} associated with $\vv$.

Call a plan $\pp$  \emph{agreeable} when $p_1 = 1$ and \emph{firm} when
$p_4 = 0$. That is, an agreeable plan always responds to $cc$ with $c$ and a firm plan always responds to $dd$ with $d$.
If both $\pp$ and $\qq$ are agreeable, then $\{ cc \}$ is a terminal set for the Markov matrix $\MM$ given by (\ref{3}) and
so  $\vv = (1,0,0,0)$  is a stationary distribution with fixation at $cc$. If both $\pp$ and $\qq$ are firm, then
$\{ dd \}$ is a terminal set for  $\MM$  and
 $\vv = (0,0,0,1)$ is a stationary distribution  with fixation at $dd$. Any convex combination of agreeable plans
(or firm plans)
is agreeable (respectively, firm).

An agreeable plan together with initial cooperation is called an \emph{agreeable strategy}.

The plans $TFT  = (1,0,1,0)$ and $Repeat  = (1,1,0,0)$ are each agreeable and firm.
The same is true for any mixture of these.  If both X and Y use $TFT$, then the outcome is determined by the initial play.
Initial outcomes $cc$ and $dd$ lead to immediate fixation. Either $cd$ or $dc$
 results in period 2 alternation between these two states. Thus,
$\{ cd,dc \}$ is another terminal set with stationary distribution $(0, \frac{1}{2},  \frac{1}{2}, 0)$. If
$a \cdot TFT + (1-a) Repeat$ is used instead
  by either player (with $0 < a < 1$ ), then eventually fixation at $cc$ or $dd$ results. There are then only two
terminal sets instead of three. The period 2 alternation described above illustrates why we needed the Cesaro limit, i.e.
the limit of averages, in (\ref{3ab}) rather than the limit per se.

Because so much work had been done on this Markov model, the
exciting new ideas in Press and Dyson [15] took people by
surprise. They have inspired a number of responses, e.g., Stewart
and Plotkin [16] and especially, Hilbe, Nowak and Sigmund [8]. I
would here like to express my gratitude to Karl Sigmund whose kind,
but firm, criticism of the initial draft directed me to this recent
work. The result is both a substantive and expository improvement.

Our purpose here is to use these new ideas to characterize the  plans that are good in the following sense.

\begin{df}\label{df1.1} A plan $\pp$ for X is called \emph{good}  if it is agreeable and if for any
general plan chosen by Y  against it and any associated limit distribution, the expected payoffs satisfy
 \begin{equation}\label{new5}
 s_Y \ \geq \ R \qquad \Longrightarrow \qquad s_Y \ = \ s_X \ = \ R.
 \end{equation}
 The plan is called \emph{of Nash type} if it is agreeable and if the expected payoffs against any Y general plan satisfy
 \begin{equation}\label{new5a}
 s_Y \ \geq \ R \qquad \Longrightarrow \qquad s_Y \ = \  R .
 \end{equation}\end{df}

A \emph{good strategy} is a good plan together with initial cooperation.
\vspace{.5cm}

By Proposition \ref{newprop1}, $s_Y = s_X = R$ iff the associated limit distribution is $(1, 0, 0, 0)$.
 In the memory-one case, $(1, 0, 0, 0)$ is a stationary distribution iff both plans are agreeable.
It is the unique stationary distribution iff, in addition, the
matrix $\MM$ is convergent. If $\pp$ is not agreeable, then
(\ref{new5}) can be vacuously true. For example, if X plays $AllD  = (0, 0, 0, 0)$, then for any Y response
$P \geq s_Y$ and the implication is true.

When both players use agreeable strategies, i.e. agreeable plans with initial cooperation, then the joint
cooperative payoff is achieved. The pair of strategies
is a Nash equilibrium exactly when the two plans are of Nash type. That is, both players receive $R$ and neither
player can do better by playing an alternative plan. A good plan is of Nash type, but
more is true. We will see that with  a Nash equilibrium it is possible that
Y can play an alternative which still yields $R$ for herself but with the payoff to X smaller than $R$. That is, Y has no incentive
to play so as to reach the joint cooperative payoff.   On the other hand, if X uses a good plan, then the only responses for Y
that obtain $R$ for her also yield $R$ for X.

The plan $Repeat = (1, 1, 0, 0)$ is an agreeable plan that is not of Nash type.
If both players use $Repeat$, then the initial outcome repeats forever. If the initial
outcome is $cd$, then $s_Y = T$ and $s_X = S$.

For a plan $\pp$, we define the \emph{X Press-Dyson vector} $\tp = \pp - \ee_{12}$, where
$\ee_{12} = (1, 1, 0, 0)$. Considering the utility
of the following result of Hilbe, Nowak and Sigmund, its proof, taken from Appendix A of [8], is remarkably simple.

\begin{theo}\label{newtheo2} Assume that X uses the plan $\pp$  with X Press-Dyson vector $\tp$. If the
initial plays and the general plan of Y yields the sequence of distributions $\{ \vv^{n} \}$, then
\begin{equation}\label{new6}
\begin{split}
Lim_{n \to \infty} \frac{1}{n} \ \Sigma_{k = 1}^n \ < \vv^k \cdot \tp > \quad = \quad 0,\hspace{2cm} \\
\mbox{and so} \quad < \vv \cdot \tp > \ = \   v_1 \tilde{p}_1 + v_2 \tilde{p}_2 + v_3 \tilde{p}_3 + v_4 \tilde{p}_4  \ = \ 0
\end{split}
\end{equation}
for any associated limit distribution $\vv$.
\end{theo}

{\bfseries Proof:} Let $v^n_{12} = v^n_1 + v^n_2$, the probability that either $cc$ or $cd$
 is the outcome in the $n^{th}$ round of play.
That is, $v^n_{12} = < \vv^n \cdot \ee_{12} > $ is the probability that X played $c$ in the $n^{th}$ round.
On the other hand, since X is
using the plan $\pp$, $p_i$ is the conditional probability that X plays $c$ in the next round, given outcome $z_i$
in the current round. Thus,
$ < \vv^n \cdot \pp >$ is the probability that X plays $c$ in the $(n + 1)^{st}$ round, i.e. it is  $v^{n+1}_{12}$.  Hence,
$v^{n+1}_{12} - v^{n}_{12} \ = \ < \vv^n \cdot \tp >$. The sum telescopes to yield
\begin{equation}\label{new7}
v^{n+1}_{12} \ - \ v^1_{12} \quad = \quad \Sigma_{k = 1}^n \ < \vv^k \cdot \tp >.
\end{equation}
As the left side has absolute value at most 1, the limit (\ref{new6}) follows. If a subsequence of the Cesaro averages converges
to $\vv$, then $< \vv \cdot \tp > \ = \ 0$ by continuity of the dot product. \hspace{4cm} $\Box$
 \vspace{.5cm}

To illustrate the use of this result, we examine $ TFT = (1, 0, 1, 0)$
 and another plan which has been labeled in the literature
$Grim  = (1, 0, 0, 0)$. We consider mixtures of each with  $Repeat  = (1, 1, 0, 0)$.

\begin{cor}\label{newcor3} Let $1 \geq a > 0$.

(a) The plan $\pp = a TFT  + (1 - a) Repeat$ is a good plan  with $s_Y = s_X$ for any limiting distribution.

(b) The plan $\pp = a Grim   + (1 - a) Repeat$ is  good.
\end{cor}

{\bfseries Proof:} (a) In this case, $\tp = a(0, -1, 1, 0)$ and so (\ref{new6}) implies that $v_2 = v_3$. Thus, $s_Y = s_X$.
From this (\ref{new5}) follows from Proposition \ref{newprop1}.

(b) Now $\tp = a(0, -1, 0, 0)$ and so (\ref{new6}) implies that $v_2 = 0$.  Thus,  $s_Y = v_1 R + v_3 S + v_4 P$ and this is
less than $ R $ unless $v_3 = v_4 = 0$ and $v_1 = 1$.  When $v_1 = 1$,  $s_Y = s_X = R$, proving (\ref{new5}). \hspace{2cm} $\Box$
\vspace{.5cm}

In the next section we will prove the following characterization of the good plans.

\begin{theo}\label{newtheo3a} Let $\pp = (p_1, p_2, p_3, p_4)$ be an agreeable plan other than $Repeat$. That is,
$p_1 = 1$ but $\pp \not= (1, 1, 0, 0)$.

The plan $\pp$ is of Nash type iff   the following inequalities hold.
\begin{equation}\label{ineq}
\frac{T - R}{R - S} \cdot p_3 \ \leq \ (1 - p_2) \qquad \mbox{and} \qquad  \frac{T - R}{R - P} \cdot p_4 \ \leq \ (1 - p_2).
\end{equation}
The plan $\pp$ is good iff, in addition, both inequalities are strict.
\end{theo}
\vspace{.5cm}

 \begin{cor}\label{newcor3b} In the compact convex set of agreeable plans,
 the set $\{ \pp $ equals $ Repeat $ or is of Nash type $ \}$ is a closed convex set with interior the set of good plans.
 \end{cor}

 {\bfseries Proof:} The X Press-Dyson vectors form a cube and the agreeable plans are the intersection
 with the subspace $\tilde{p}_1 = 0$. We then intersect with the half-spaces defined by
\begin{equation}\label{ineq2}
\frac{T - R}{R - S}  \tilde{p}_3  + \tilde{p}_2 \ \leq \ 0 \qquad \mbox{and}
\qquad  \frac{T - R}{R - P} \tilde{p}_4  + \tilde{p}_2\ \leq \ 0.
\end{equation}
The result is a closed convex set with interior given by the strict inequalities. Notice that these conditions are preserved
by multiplication by a positive constant $a \leq 1$ or by any larger constant so long as $a \tp$ remains in the
cube. Hence, $Repeat$ with $\tp = 0$ is on the boundary. \hspace{1.5cm} $\Box$
\vspace{.5cm}

It is easy to compute that that
\begin{equation}\label{new8}
det \begin{pmatrix} R & R & 1 & 0 \\
S & T & 1 & 1 \\
T & S & 1 & 1  \\
P & P & 1 & 0 \end{pmatrix} \quad = \quad - 2 (R - P)(T - S).
\end{equation}
Hence, with $\ee_{23} \ = \ (0, 1, 1, 0)$,  we can use $\{ \SX, \SY, \one, \ee_{23} \}$ as a basis for $\R^4$. For a
distribution vector $\vv$ we will write $v_{23} $ for $v_2 + v_3 \ = \ < \vv \cdot \ee_{23} >$. From Theorem \ref{newtheo2},
we immediately obtain the following.

\begin{theo}\label{newtheo4} If $\pp$ is a plan whose X Press-Dyson vector $\tp = \a \SX + \b \SY + \g \one + \d \ee_{23}$
and $\vv$ is a limit distribution when Y plays some general plan against $\pp$, then the average payoffs
satisfy the following \emph{Press-Dyson Equation}.
\begin{equation}\label{new9}
\a s_X \ + \ \b s_Y \ + \ \g \ + \ \d v_{23} \quad = \quad 0.
\end{equation}
\end{theo}
 \vspace{.5cm}

 The most convenient cases to study occur when $\d = 0$. Press and Dyson called such a plan a \emph{Zero-Determinant Strategy}
 (hereafter ZDS) because of an ingenious determinant argument leading to (\ref{new9}). We have used
 Theorem \ref{newtheo2} of  Hilbe-Nowak-Sigmund   instead.

 This representation yields a simple description of the good plans.

\begin{theo}\label{newtheo5} Assume that $\pp = (p_1, p_2, p_3, p_4)$ is an agreeable plan
 with X Press-Dyson vector
$\tp = \a \SX + \b \SY + \g \one + \d \ee_{23}$. Assume that $\pp$ is not $Repeat$, i.e. $(\a, \b, \g , \d) \not= (0, 0, 0, 0)$.
The plan $\pp$ is of Nash type iff
\begin{equation}\label{ineq}
max( \frac{ \d }{(T - S)}, \frac{ \d }{(2R - (T + S))} ) \quad \leq \quad  \a. \hspace{2cm}
\end{equation}
The plan $\pp$ is good iff, in addition, the inequality is strict.
\end{theo}
\vspace{.5cm}

{\bfseries Remark:} Observe that $ T - S \ > \ 2R - (T + S) \ > \ 0$.   It follows that if $\d \leq 0$, then
$\pp$ is good iff  $\frac{ \d }{(T - S)}  < \a$. On the other hand,  if $\d > 0$, then $\pp $ is good iff
$\frac{ \d }{(2R - (T + S))} < \a$.
\vspace{.5cm}

In the next section, we will investigate the geometry of the $\{ \SX, \SY, \one, \ee_{23} \}$ decomposition
of the Press-Dyson vectors and prove the theorems.
\vspace{1cm}

\section{ Good Plans and   The \\ Press-Dyson  Decomposition }

We begin by normalizing the payoffs. We can add to all a common number and
 multiply all by a common positive number without changing the relationship between the various strategies. We subtract
 $S$ and divide by $T - S$. So from now on we will assume
 that $T = 1$ and $S = 0$.

The payoff vectors of (\ref{1a}) are then given by
\begin{equation}\label{14}
\SX \quad = \quad (R, 0, 1, P), \qquad \SY \quad = \quad (R, 1, 0, P),
\end{equation}
and from (\ref{2}) we have
\begin{equation}\label{13}
 1 \ > \ R \ > \  \frac{1}{2}, \quad \mbox{and} \quad R \ > \ P \ > \ 0.
\end{equation}

After normalization Theorem \ref{newtheo3a} becomes the following.

\begin{theo}\label{newtheo3anorm} Let $\pp = (p_1, p_2, p_3, p_4)$ be an agreeable plan other than $Repeat$.
That is, $p_1 = 1$ but $\pp \not= (1, 1, 0, 0)$.

The plan $\pp$ is of Nash type iff   the following inequalities hold.
\begin{equation}\label{ineq2}
\frac{1 - R}{R} \cdot p_3 \ \leq \ (1 - p_2) \qquad \mbox{and} \qquad  \frac{1 - R}{R - P} \cdot p_4 \ \leq \ (1 - p_2).
\end{equation}
The plan $\pp$ is good iff, in addition, both inequalities are strict.
\end{theo}
\vspace{.5cm}

{\bfseries Proof:} We first eliminate the possibility $p_2 = 1$. If $1 - p_2 = 0$, then the inequalities would yield
$p_3 = p_4 = 0$ and so $\pp = Repeat$, which we have excluded.  On the other hand, if $p_2 = 1$, then
$\pp = (1, 1, p_3, p_4)$. If against this Y plays $AllD  = (0, 0, 0, 0)$, then $\{ cd \}$ is a terminal set with stationary
distribution $(0, 1, 0, 0)$ and so with $s_Y = 1$ and $s_X = 0$. Hence, $\pp$ is not of Nash type.  Thus, if
$p_2 = 1$, then neither is $\pp$  of Nash type, nor do the inequalities hold for it.  We now assume $1 - p_2 > 0$.

Observe that
\begin{equation}\label{proof1}
\begin{split}
s_Y  -  R   \quad = \quad (v_1 R + v_2 + v_4 P) \ - \ (v_1 R + v_2 R + v_3 R + v_4 R) \\
= \quad v_2 (1 - R) \ - \ v_3 R \ - \ v_4 (R - P). \hspace{2cm}
\end{split}
\end{equation}
Hence, multiplying by the positive quantity $(1 - p_2)$, we have
\begin{equation}\label{proof2}
s_Y \ >= \ R \quad \ \Longleftrightarrow \quad \ (1 - p_2) v_2 (1 - R) \ >= \ v_3 (1 - p_2)R  \ + \ v_4 (1 - p_2) (R - P),
\end{equation}
where  this notation means that the inequalities are equivalent and the equations are equivalent.

Since $\tilde{p}_1 = 0$, equation (\ref{new6}) implies $v_2 \tilde{p}_2 + v_3 \tilde{p}_3 +  v_4 \tilde{p}_4  = 0$
and so $(1 - p_2) v_2 = v_3 p_3 + v_4 p_4$. Substituting in the above inequality and collecting terms we get
\begin{equation}\label{proof3}
\begin{split}
 s_Y \ >= \ R \qquad \Longleftrightarrow \qquad   A v_3 \  > =  \ B v_4 \qquad \mbox{with} \hspace{3cm} \\
A \ = \ [ p_3 (1 - R) - (1 - p_2) R] \qquad \mbox{and} \qquad B \ = \ [(1 - p_2)(R - P) - p_4 (1 - R)].
\end{split}\end{equation}
Observe that the inequalities of (\ref{ineq2}) are equivalent to $A \leq 0$ and $B \geq 0$.
The proof is completed by using a sequence of little cases.
\vspace{.25cm}

Case(i) $ A = 0, B = 0$ : In this case, $A v_3 =  B v_4$ holds for any strategy for Y. So for any Y strategy, $s_Y = R$ and
$\pp$ is of Nash type.  If Y chooses a plan that is not agreeable, then $\{ cc \}$ is not a closed set of states and
so $v_1 \not= 1$.  From Proposition \ref{newprop1}, $s_X < R$ and so $\pp$ is not good.
\vspace{.25cm}

Case(ii) $ A < 0, B = 0$ : The inequality $A v_3 \geq  B v_4$ holds
iff $v_3 = 0$. If $v_3 = 0$, then  $A v_3 =  B v_4$ and so $s_Y = R$. Thus, $\pp$ is  Nash.
\vspace{.25cm}

Case(iia) $B \leq 0$, any $A$: Assume Y chooses a plan that is not agreeable and is such  that $v_3 = 0$.   For example,
if Y plays $AllD = (0, 0, 0, 0)$, then no state moves to $dc$. With such a Y choice,  $A v_3 \geq  B v_4$ and so $s_Y \geq R$.
As above,  $v_1 \not=1$ because the Y plan is not agreeable. Again $s_X < R$ and $\pp$ is not good.
Furthermore, $v_3 = 0, v_1 < 1, p_2 < 1,$ and $ (1 - p_2) v_2 = v_4 p_4$ imply
that $v_4 > 0$. So if $B < 0$, then  $A v_3 >  B v_4$ and so $s_Y > R$. Hence, $\pp$ is not Nash when $B < 0$.
\vspace{.25cm}

Case(iii) $A = 0, B > 0$: The inequality $A v_3 \geq  B v_4$ holds iff $v_4 = 0$. If
 $v_4 = 0$, then  $A v_3 =  B v_4$ and $s_Y = R$. Thus, $\pp$ is Nash.
\vspace{.25cm}

Case(iiia) $A \geq 0$, any $B$: Assume Y  chooses a plan that is not agreeable and is such that
 $v_4 = 0$. For example, if Y plays $(0, 1, 1, 1)$,  then no state moves to $dd$. With such a Y choice,
  $A v_3 \geq  B v_4$ and so $s_Y \geq R$. As before,
  $v_1 \not=1 $ implies $s_X < R$ and the plan is not good.  Furthermore,
  $v_4 = 0, v_1 < 1, p_2 < 1,$ and $ (1 - p_2)v_2 = v_3 p_3$  imply
  that $v_3 > 0$. So if $A > 0$, then  $A v_3 >  B v_4$ and so $s_Y > R$. Hence, $\pp$ is not Nash when $A > 0$.
\vspace{.25cm}

 Case(iv) $A < 0, B > 0$:  The inequality $A v_3 \geq  B v_4$ implies $v_3, v_4 = 0$. So $(1 - p_2) v_2 = v_3 p_3 + v_4 p_4 = 0$.
 Since $p_2 < 1, \ v_2 = 0$.
 Hence, $v_1 = 1$.  That is, $s_Y \geq R$ implies $s_Y = s_X = R$ and so $\pp$ is good. \hspace{2cm}  $\Box$
 \vspace{.5cm}

{\bfseries Remarks:} (a) Since $1 > R > \frac{1}{2}$, it is always true that $\frac{1 - R}{R} < 1$. On the other hand,
$\frac{1 - R}{R - P}$ can be greater than $1$ and  the second inequality requires $p_4 \leq \frac{R - P}{1 - R}$.
In particular, if $p_2 = 0$, then the plan is good iff  $p_4 < \frac{R - P}{1 - R}$.
For example, the plan $(1, 0, 0, 1)$ is,  in the literature, labeled $Pavlov$, or $Win Stay, Lose Shift$.
This plan always satisfies the
first inequality strictly, but it satisfies the second strictly, and so is good, iff $1 - R < R - P$.

(b) In Case(i) of the proof, the payoff $s_Y = R$ is determined by
$\pp$ independent of the choice of strategy for $Y$.  In general, plans
 that fix the opponent's payoff in this way were described
 by Press and Dyson [15] and, earlier, by Boerlijst, Nowak and Sigmund [5], where they are called  \emph{equalizer strategies}.
 The agreeable equalizer plans have $\tp = a (0, - \frac{1-R}{R}, 1, \frac{R - P}{R})$ with $1 \geq a > 0$.
\vspace{.5cm}

Christian Hilbe suggests a nice interpretation of the above results:

\begin{cor}\label{Hilbecor} Let $\pp$ be an agreeable plan with $p_2 < 1$.

(a) If $\pp$ is good, then using any plan $\qq$ that is not agreeable forces Y to get a payoff $s_Y < R$.

(b) If $\pp$ is not good, then by using at least one of the two plans $\qq = (0, 0, 0, 0)$ or $\qq = (0, 1, 1, 1)$, Y
can certainly obtain a payoff $s_Y \geq R$, and force X to get a payoff $s_X < R$.

(c) If $\pp$ is not Nash, then by using at least one of the two plans $\qq = (0, 0, 0, 0)$ or $\qq = (0, 1, 1, 1)$ Y
can certainly obtain a payoff $s_Y > R$, and force X to get a payoff $s_X < R$.
\end{cor}

{\bfseries Proof:} (a) If $\pp$ is good, then $s_Y \geq R$ implies $s_Y = s_X = R$ which requires $v = (1, 0, 0, 0)$.
This is only stationary when $\qq$ as well as $\pp$ is agreeable.

(b) and (c) follow from the analysis of  cases in the above proof. \hspace{1cm}  $\Box$
\vspace{.5cm}

{\bfseries Remark:} If $p_2 = p_1 = 1$, then the plan $\pp$ is not Nash. As observed in the proof of Theorem \ref{newtheo3anorm}
above, if
Y plays $\qq = (0, 0, 0, 0)$,
then $cd$ is a terminal set with stationary distribution $\vv = (0, 1, 0, 0)$ and so with $s_Y = 1, s_X = 0$.
However, if, in addition, $p_4 = 0$, e.g. if X uses $Repeat$, then $dd$ is also a terminal set. Thus,
 if X plays $\pp$ with $1 - p_4 = p_2 = p_1 = 1$ and Y always defects,
then fixation occurs immediately at either $cd$ with $s_Y = 1$ and $s_X = 0$, or else at $dd$ with $s_Y = s_X = P$.
The result is  determined by the initial play of X.
\vspace{.5cm}

We now consider the Press-Dyson representation, using the normalized payoff vectors of (\ref{14}).
If $\tp = \a \SX + \b \SY + \g \one + \d \ee_{23}$ is the X Press-Dyson vector of a plan $\pp$, then
it must satisfy two sorts of constraints.

The \emph{sign constraints} require that the first two entries be nonpositive and the last two be nonnegative. That is,
\begin{equation}\label{15}
\begin{split}
(\a + \b) R \ + \ \g  \quad \leq \quad 0, \\
\qquad \b \ + \ \g  \ + \ \d \quad \leq \quad 0,\\
\qquad \a \ + \ \g   \ + \ \d \quad \geq \quad 0,\\
(\a + \b) P \ + \ \g \quad \geq \quad 0.
\end{split}
\end{equation}

\begin{lem}\label{newlem2.1} If $\tp = \a \SX + \b \SY + \g \one + \d \ee_{23}$ satisfies the sign
constraints, then
\begin{equation}\label{15new}
\begin{split}
\a  + \b \ \leq \ 0  \qquad \mbox{and} \qquad \g \geq 0, \\
\a  + \b \ = \ 0  \qquad \Leftrightarrow \qquad \g = 0.
\end{split}
\end{equation}
\end{lem}

{\bfseries Proof:} Subtracting the fourth inequality from the first we see that $(\a + \b)(R - P) \leq 0$ and so
$R - P > 0$ implies $\a + \b \leq 0$. Then the fourth inequality and $P > 0$ imply $\g \geq 0$.
The first and fourth imply $\a + \b = 0 $ iff $\g = 0$. \hspace{2cm} $\Box$
\vspace{.5cm}

{\bfseries Remark:} Notice that both $\tilde{p}_1$ and $ \tilde{p}_4 $ vanish iff $\a + \b = \g = 0$.
These are the cases when the plan $\pp$ is both agreeable and firm.
\vspace{.5cm}

In addition, the entries of an X Press-Dyson vector have absolute value at most 1.
These are the \emph{size constraints}. If a vector satisfies the sign constraints then,
 multiplying by a sufficiently small positive number, we obtain
the size constraints as well.
  Any vector in $\R^4$ that satisfies both the sign and the size constraints
 is an X  Press-Dyson vector. Call $\pp$ a \emph{top plan}
 if $|\tilde p_i| = 1$ for some $i$. For
 any plan $\pp$, other than $Repeat$, which has X Press-Dyson vector ${\mathbf 0}$,
 $\pp = a(\pp^t) + (1 - a)Repeat$ for a unique top plan $\pp^t$ and a unique positive $a \leq 1$.
 Equivalently, $\tp = a \tp^t$.

 Observe that $\pp$ is agreeable iff  $\tilde{p}_1 = 0$ and so iff $(\a + \b)R + \g = 0$. In that case,
   $\b = - \a - \g R^{-1}$.
 Substituting into (\ref{new9}), we obtain the following corollary of Theorem \ref{newtheo4}.

 \begin{cor}\label{newcor2.2}   If $\pp$ is an agreeable plan with
 X Press-Dyson vector $\tp = \a \SX + \b \SY + \g \one + \d \ee_{23}$,
 then the payoffs
with any limit distribution
satisfy the following version of the \emph{Press-Dyson Equation}.
 \begin{equation}\label{new16}
 \g R^{-1} s_Y \ + \ \a (s_Y - s_X) \ - \ \d v_{23} \quad = \quad \g.
 \end{equation}
 \end{cor}
  \vspace{.5cm}

  Now we justify the description in Theorem \ref{newtheo5}. Notice that if we label
  by $\SX^0$ and $\SY^0$ our original payoff vectors before normalization then
  $\SX^0 = (T - S) \SX + S \one, \SY^0 = (T - S) \SY + S \one$ and so if $(\a,\b,\g,\d)$
  are the coordinates of $\tp$ with respect to the basis $\{ \SX, \SY, \one, \ee_{23} \}$
  then $(\a^0,\b^0,\g^0,\d^0) = (\a/(T-S),\b/(T-S), \g - (\a + \b)S/(T-S), \d)$ are the
  coordinates with respect to $\{ \SX^0, \SY^0, \one, \ee_{23} \}$. In particular,
$ \a >= k \d$ iff $\a^0 >= k \d^0/(T-S)$  for any $k$. Furthermore, the constant $k = \frac{T - S}{2R - (T + S)}$ is
independent of normalization. So it suffices to prove the normalized version of the theorem which is the following.

\begin{theo}\label{newtheo5norm} Assume that $\pp = (p_1, p_2, p_3, p_4)$ is an agreeable plan
 with X Press-Dyson vector
$\tp = \a \SX + \b \SY + \g \one + \d \ee_{23}$. Assume that $\pp$ is not $Repeat$,
i.e. $(\a, \b, \g , \d) \not= (0, 0, 0, 0)$.
 The plan $\pp$ is of Nash type iff
\begin{equation}\label{ineqnorm}
max(  \d, (2R - 1)^{-1} \d ) \quad \leq \quad  \a.   \hspace{4cm}
\end{equation}
The plan $\pp$ is good iff, in addition, the inequality is strict.
\end{theo}
\vspace{.5cm}

{\bfseries Proof:} Since  $\b = - \a - \g R^{-1}$, we have
\begin{equation}
\begin{split}
(1 - p_2) \ = \ - \tilde{p}_2 \ = \ - \b - \g - \d \ = \ \a + \frac{1 - R}{R}\g - \d, \\
p_3 \ = \ \tilde{p}_3 \ = \ \a + \g + \d, \qquad p_4 \ = \ \tilde{p}_4 \ = \ \frac{R - P}{R}\g.
\end{split}
\end{equation}

The inequality $(1 - R)p_3 \leq R (1 - p_2)$ becomes $(1 - R)(\a + \g + \d) \leq R \a + (1 - R) \g -R \d$.
 This reduces to $\d \leq (2R - 1) \a$. Similarly, the inequality $(1 - R) p_4 \leq (R - P) (1 - p_2) $
reduces to $ \d  \leq \a$. \hspace{1.5cm} $\Box$
\vspace{.5cm}

{\bfseries Remarks:} (a) Thus, when $\d \leq 0$, $\pp$ is good iff $\d < \a$. When $\d > 0$, $\pp$ is
good iff $\frac{ \d }{2R - 1} < \a$.

(b) From the proof we see that the equalizer case, when both inequalities of (\ref{ineq2}) are equations,
occurs when $\d = \a = (2R - 1)^{-1} \d$.  Since $2R - 1 < 1$, this reduces to $0 = \d = \a$.
\vspace{.5cm}

In the ZDS case, when $\d = 0$,  we can rewrite (\ref{new16}) as
\begin{equation}\label{new16a}
\kappa \cdot ( s_X - R) \quad = \quad   s_Y - R  \hspace{3cm}
\end{equation}
with
$\kappa = \frac{\a R}{\g + \a R}$. Thus, the condition $\a >  0$ is equivalent to $0 < \kappa \leq 1$. In [8]
these plans are introduced and called \emph{complier strategies}. The equation and the condition $\kappa > 0$ make
it clear that such plans are good. In addition, if $s_Y < R$, then it follows that $s_X \leq s_Y$ with strict
inequality when $\g > 0$ and so $\kappa < 1$. The strategy ZGTFT-2 analyzed in Stewart and Plotkin [16] is an
example of a complier plan.
 When X plays a complier plan, then either both $s_X$ and $s_Y$ are equal to $R$, or
else both are below $R$.  This is not true for good plans in general.  If X plays the good plan $Grim = (1, 0, 0, 0)$
and Y plays $(0, 1, 1, 1)$, then fixation at $dc$ occurs with $v = (0, 0, 1, 0)$ and so with $s_Y = 0 $ ($< R$ as required by
Corollary \ref{Hilbecor} ), but with  $s_X = 1 > R$.
\vspace{.5cm}

 Let us look at the geometry of the Press-Dyson representation.

We begin with the \emph{exceptional plans} which are defined by $\g = \a + \b = 0$. The sign constraints
yield $\a = -  \b \geq | \d |$ and $\tp = (0, \d - \a, \d + \a, 0)$. As remarked after Lemma \ref{newlem2.1},
the exceptional plans are exactly those plans that
are both agreeable and firm. In the $xy$ plane with $x = \tilde{p}_2$ and $y = \tilde{p}_3$ they form a square with
vertices:

 $Repeat$ ( $\tp = (0,0,0,0)$ ), $Grim$ ( $\tp = (0,-1,0,0)$ ),

 $TFT$ ( $\tp = (0,-1,1,0)$ ), and
what we will call $Lame$ ( $\tp = (0,0,1,0)$). \\
Thus, $Lame  = (1,1,1,0)$. The top plans consist of the segment
that connects  $TFT$ with $Grim$ together with the segment that connects  $TFT$ with $Lame$.

On the $Grim-TFT$ segment $\d - \a = -1$ and $0 \leq \d + \a \leq 1$. That is, $\d = \a - 1$ and $-\frac{1}{2} \leq \d \leq 0$.
By Theorem \ref{newtheo5norm}, the plans in the triangle with vertices Grim, TFT and Repeat are all good except
for Repeat itself.

On the $Lame-TFT$ segment $-1 \leq \d - \a \leq 0$ and $ \d + \a = 1$. That is, $\d = 1 - \a$ and $\frac{1}{2} \geq \d \geq 0$.
Such a plan is the mixture $t TFT + (1 -t)Lame = (1, 1 - t, 1, 0)$ with $t = 2 \a - 1$. By Theorem \ref{newtheo3anorm},
this plan is good  iff $t > \frac{1 - R}{R} $.
The plan on the $TFT-Lame$ segment with $t = \frac{1 - R}{R} $, and so with $2 \a = R^{-1}$, we will call
$Edge = (1,\frac{2R - 1}{R},1,0)$. The plans in the $TFT-Edge-Repeat$ triangle that are not on
the $Edge-Repeat$ side are good plans. The  plans in the complementary
$Edge-Lame-Repeat$ triangle are not good.
\vspace{.5cm}

Now assume $\g > 0$, and   define
\begin{equation}\label{15a}
\bar \a  \quad = \quad \a / \g, \qquad \bar \b \quad = \quad \b/ \g, \qquad \bar \d \quad = \quad \d/ \g.
\end{equation}
with the sign constraints
\begin{equation}\label{16a}
\begin{split}
-P^{-1} \quad \leq \quad \bar \a  \ + \ \bar \b \quad \leq \quad - R^{-1},\hspace{1cm} \\
\bar \b  \quad \leq \quad -1 - \bar \d \quad \leq \quad \bar \a. \hspace{2cm}
\end{split}
\end{equation}
For any triple $(\bar \a,\bar \b, \bar \d)$ that satisfies these
inequalities, we obtain an X Press-Dyson vector vector $\tp = \a \SX
+ \b \SY + \g \one + \d \ee_{23}$ that satisfies the size
constraints as well by using $(\a, \b, \g, \d) \ = \ \g \cdot (\bar
\a,\bar \b, 1, \bar \d)$ with $\g
> 0$ small enough. When we use the largest value of $\g$ such that
the size constraints hold, we obtain the top plan associated
with $(\bar \a,\bar \b, \bar \d)$. The others are mixtures of the
top plan with Repeat. For a plan with this triple the
Press-Dyson Equation (\ref{new9}) becomes
\begin{equation}\label{newpressdyson}
\bar \a s_X \ + \ \bar \b s_Y \ + \ \bar \d v_{23} \ + \ 1 \quad = \quad 0.
\end{equation}

The points $(x,y) = (\bar \a,\bar \b)$ lie in the \emph{Strategy Strip}.  This consists of the points of the $xy$ plane
with $y \leq x$ and
that lie on or below the line $x + y = -R^{-1}$ and on or above the line $x + y = -P^{-1}$. Then
$\bar \d$ must satisfy $-1 -x \leq \bar \d \leq -1 -y$. Alternatively, we can fix $\bar \d$ to be arbitrary and
intersect the Strategy Strip with the fourth quadrant when the origin is at $(-1 - \bar \d, -1 - \bar \d)$, i.e. the
points with $y \leq -1 - \bar \d \leq x$.

Together with the exceptional plans those with $(\bar \a, \bar \b)$  on the line $x + y = -R^{-1}$ are
exactly the agreeable plans.
Together with the exceptional plans those on the line
$x + y = -P^{-1}$ are exactly the firm plans.
\vspace{.5cm}

Let us look at the good ZDS's, i.e. the good plans with $\d = 0$. In the exceptional case with $\g = 0 $, the
top good plan is $TFT$. When $\d = 0$ and $\g > 0$,  the good plans are those that satisfy
$\bar \a + \bar \b = - R^{-1} $ and $\bar \a > 0$.  As mentioned above, these are the
complier plans.

\begin{prop}\label{prop2.4b} Given $\bar \a > 0$, the associated agreeable ZDS top plan is given by
\begin{equation}\label{22}
 \pp \quad = \quad (1,  \frac{2R  - 1}{R(\bar \a + 1)}, 1, \frac{R - P }{R(\bar \a + 1)}).
\end{equation}
\end{prop}

{\bfseries Proof:}  The agreeable plan $\pp$ with $\g, \bar \a > 0$ and $\bar \d = 0$ has X Press-Dyson vector
\begin{equation}\label{22aa}
\tp \quad = \quad (0, -\g(\bar \a + R^{-1} - 1), \g(\bar \a + 1), \g( 1 - P \cdot R^{-1})).
\end{equation}
  With $\bar \a$ fixed, the largest
 value for $\g$ so that the size constraints hold is $(\bar \a + 1)^{-1}$ .
 This easily yields (\ref{22}) for the top plan.\hspace{1cm} $\Box$
 \vspace{.5cm}

When $\bar \d = 0$, the vertical line $\bar \a = 0$ intersects the strip in points
whose plans are all the \emph{equalizers}, as
discussed by Press and Dyson [15] and by Boerlijst, Nowak and Sigmund [5].  Observe that with $\bar \d = 0$ and
$\bar \a = 0$ the Press-Dyson Equation (\ref{newpressdyson}) becomes
$ \bar \b s_Y \ +  \ 1 \  = \  0, $ and so $s_Y  =  - \bar \b^{-1}$ regardless of the choice of
strategy for Y. The agreeable case has $\bar \b = -R^{-1}$.
 The vertical line of equalizers cuts the line of agreeable plans, separating it into the unbounded
ray with good plans and the segment with plans that are not even of Nash type.

Finally, we call a plan $\pp$ \emph{generous} when $p_2 > 0$ and $p_4 > 0$. That is, whenever Y
defects there is a positive probability that X will cooperate. The complier plans given by (\ref{22})
are generous.

\begin{prop}\label{prop2.6} Assume that X plays $\pp$, a generous plan of Nash type. If Y plays
plan $\qq$ of Nash type and either (i) $\qq$ is generous, or (ii) $q_3 + q_4 > 0$, then
 $\{cc \}$ is the unique terminal set for the associated Markov matrix $\MM$. Thus, $\MM$ is convergent.
\end{prop}

{\bfseries Proof:}  Since $\pp$ and $\qq$ are both agreeable, $\{ cc \}$ is a terminal set for $\MM$.

Since $\pp$ is of Nash type, it is not $Repeat$ and so (\ref{ineq2}) implies that $p_2 < 1$.

For the first case, we prove that if $p_1 = 1, p_2 < 1, p_4 > 0$ and $\qq$ satisfies analogous conditions
and not both $\pp$ and $\qq$ are of the form $(1, 0, 1, a)$, then $\MM$ is convergent.

Recall that Y responds to $cd$ using $q_3$ and to $dc$ using $q_2$.

The assumptions $p_4, q_4 > 0$ imply that there is an edge from $dd$ to $cc$, and so that $dd$ is transient. There is an edge from
$dc$ to $dd$ if $p_3 < 1$ since $q_2 < 1$. If $p_3 = 1$ and $q_2 > 0$, then there is an edge to $cc$. There remains
the case that $p_3 = 1, q_2 = 0$ with the only edge from $dc$ going to $cd$. Similarly, there is an
edge from $cd$ to either $dd$ or $cc$ except when $p_2 = 0, q_3 = 1$. Thus, the only case when $\MM$ is not convergent
is when both $\pp$ and $\qq$ are of the form $(1, 0, 1, a)$. In that case, $\{ cd, dc \}$ is an additional terminal set.
In particular, if either $p_2$ or $q_2$ is positive, then $\{ cc \}$ is the only terminal set.  This completes
case (i). It also shows that if $\pp$ is generous and $q_4 >0$, then $\MM$ is convergent.

To complete case(ii), we assume that $\pp$ is generous and $q_3 > 0$. Since $p_2 > 0$ and $q_3 > 0$, there is an
edge from $cd$ to $cc$ and so $cd$ is transient. Since $p_4 > 0$, there is an edge from $dd$ either to $cc$ or to
$cd$ and so $dd$ is transient. Finally, $q_2 < 1$ implies there is an edge from $dc$ to $cd$ or to $dd$.  Thus, $dc$ is
transient as well. \hspace{4cm} $\Box$
 \vspace{.5cm}

This result indicates the advantage which the good plans that are generous have over the good exceptional plans like
$Grim$ and $TFT$. The latter are firm as well as agreeable.  Playing them against each other yields a nonconvergent matrix
with both $\{ cc \}$ and $\{ dd \}$ as terminal sets.  Initial cooperation does lead to immediate fixation at
$cc$, but an error might move  the sequence of outcomes on a path leading to another terminal set. When generous good plans
are used against each other, $\{ cc \}$ is the unique terminal set.  Eventual fixation at $cc$ occurs whatever the initial
distribution is, and if an error occurs, then the strategies move the successive outcomes along a path that returns to $cc$.
 It is easy to compute
 the expected number of steps $T_z$ from transient state $z$ to $cc$.
 \begin{equation}\label{32aa}
 T_z \quad = \quad 1 \ + \ \Sigma_{z'} p_{zz'} T_{z'},
 \end{equation}
where we sum over the three transient states and $p_{zz'}$ is the probability of moving along an edge from $z$ to $z'$.
Thus, with $\MM' = \MM -I$, we obtain the formula for the vector ${\mathbf T} = (T_2,T_3,T_4)$:
\begin{equation}\label{33aa}
\MM'_t \cdot {\mathbf T} \quad =\quad - {\mathbf 1}.
\end{equation}
where $\MM'_t$ is the invertible $3 \times 3$ matrix obtained from $\MM'$ by omitting the first row and column.

Consider the case when X and Y both use the  plan  given by (\ref{22}), so that $\pp = \qq = (1,p_2,1,p_4)$.
The only edges coming from
$cd$ connect with $cc$ or with $dc$ and similarly for the edges from $dc$.  Symmetry will imply that
 $T_{cd} = T_{dc}$. So with $T$ this common value we obtain from (\ref{32aa})
$T = 1 + (1 - p_2)T$.   Hence, from (\ref{22}) we get
\begin{equation}\label{34aa}
T \ = \ T_{cd} \ = \ T_{dc} \quad = \quad \frac{1}{p_2} \ = \ \frac{\bar \a + 1}{2 - R^{-1}}.
\end{equation}
Thus, the closer the plan is to the equalizer plan with $ \bar \a = 0$ the shorter the expected
recovery time from an error leading to a $dc$ or $cd$ outcome. From (\ref{32aa}) one can see that
\begin{equation}\label{35aa}
T_{dd} \quad = \quad 1 \ + \ 2p_4(1 - p_4)\cdot T \ + \ (1 - p_4)^2 \cdot T_{dd}.
\end{equation}
We won't examine this further as arriving at $dd$ from $cc$ implies errors on the part of both players.

Of course, one might regard such departures from cooperation not as noise or error but as ploys.  Y might try a rare move to
$cd$ in order to pick up the temptation payoff for defection as an occasional bonus.  But if this is strategy
rather than error, it means that Y is departing from the good  plan to one with $q_1$ a bit less than $1$. Corollary
\ref{Hilbecor}(a) implies that Y loses by executing such a ploy.
\vspace{1cm}

\section{Competing Zero Determinant Strategies}

We now examine the ZDS's in more detail. Recall that a plan $\pp$ is a ZDS when $\d = 0$ in the
Press-Dyson decomposition of the X Press-Dyson vector $\tp = \pp -\ee_{12}$. With the normalization (\ref{14})
the inverse matrix of
$( \ \SX \ \SY \ \one \ \ee_{23} \ )$ is
\begin{equation}\label{invmatrix}
\frac{-1}{2(R - P)} \
 \begin{pmatrix} -1 & R - P & P - R  & 1 \\
-1 & P - R  &  R - P & 1 \\
2P & 0 & 0 & -2R  \\
 1 - 2P & P - R  & P - R  & 2R - 1 \end{pmatrix}
\end{equation}
and so if $\tp =  \a \SX \ + \  \b \SY \ + \ \g \one \ + \ \d \ee_{23}$,
\begin{equation}\label{delta}
2(R - P) \d \ = \ ( 2P - 1)\tilde p_1 \ + \ (R - P) (\tilde p_2 + \tilde p_3) \ - \ (2R - 1)\tilde p_4.
\end{equation}
Thus, for example, if $R + P = 1$,  both $AllD$ with $\tp = (-1, -1, 0, 0)$ and $AllC$ with $\tp = (0, 0, 1, 1)$ are ZDS.

The exceptional ZDS's, which have
$\g = 0$ as well as $\d = 0$, are mixtures of $TFT$ and $Repeat$.  Otherwise, $\g > 0$ and we can write
$\tp = \g(\bar \a \SX \ + \ \bar \b \SY \ + \ \one)$. When $(\bar \a, \bar \b)$ lies in the \emph{ZDSstrip} defined by
\begin{equation}\label{zdsstrip}
ZDSstrip \quad = \quad \{  (x,y) : x \geq -1 \geq y \  \mbox{and} \  -R^{-1} \geq x + y \geq -P^{-1} \ \},
\end{equation}
then the sign constraints are satisfied. The size constraints hold as well when $\g > 0$ is small enough. For $Z$
with $P \leq Z \leq R$ the intersection of the ZDSstrip with the line $x + y = -Z^{-1}$ is a \emph{value line} in the
strip.

\begin{lem}\label{lem3.1} Assume that $(\bar \a, \bar \b)$ in the ZDS strip, with $\bar \a +  \bar \b = -Z^{-1}$.  We then have
$ - \bar \b \geq  max(1,| \bar \a |)$ and  $ - \bar \b =  | \bar \a |$
iff $\bar \a = \bar \b = -1$.  If $(\bar a,\bar b)$ is also in the strip then $D = \bar \b \bar b \ - \ \bar \a \bar a \geq 0$
with equality iff $\bar \a = \bar \b = \bar a = \bar b = -1$.\end{lem}

{\bfseries Proof:}  By definition of $Z$,  $- \bar \b = \bar \a + Z^{-1} > \bar \a$.
Also, the sign constraints imply $- \bar \b \geq 1 \geq - \bar \a$,  and so $- \bar \b  \geq - \bar \a $
with equality iff $\bar \a = \bar \b = -1$.
$D \geq (- \bar \b)(- \bar b) \ - \ |\bar \a | |\bar a | \geq 0$ and the inequality is strict unless
$\bar \a = \bar \b = \bar a = \bar b = -1$.   \hspace{3cm} $\Box$
\vspace{.5cm}

{\bfseries Remark:} Because $R > \frac{1}{2}$ it is always true that $-R^{-1} > -2 $, but $-2 \geq -P^{-1}$ iff
$\frac{1}{2} \geq P$. Hence, $(-1,-1)$ is in the ZDSstrip iff $\frac{1}{2} \geq P$.
\vspace{.5cm}

For a ZDS we can usefully transform  the Press-Dyson equation (\ref{newpressdyson}).

\begin{prop} Assume that X uses  plan $\pp$ with X Press-Dyson vector
$\tp = \g(\bar \a \SX \ + \ \bar \b \SY \ + \ \one), \ \g > 0$. Let $-Z^{-1} = \bar \a + \bar \b$ , so that
 $P \leq Z \leq R$.

For any general plan played by Y,
\begin{equation}\label{newnewpressdyson2}
 \bar \a Z (s_X \ - \  s_Y) \quad = \quad (s_Y \ - \ Z).
\end{equation}

If $\kappa =  \bar \a Z/(1 + \bar \a Z)$, then $1 > \kappa $ and $\kappa$ has the same sign as $\bar \a$.
For any general plan played by Y,
\begin{equation}\label{newnewpressdyson}
 \kappa (s_X \ - \ Z) \quad = \quad (s_Y \ - \ Z).
\end{equation}
\end{prop}

{\bfseries Proof:} Notice that $1 + \bar \a Z = - \bar \b Z \geq Z  \geq P > 0$. Multiplying (\ref{newpressdyson}) by
$Z$ and substituting for $\bar \b Z$ easily yields (\ref{newnewpressdyson2}) and then (\ref{newnewpressdyson}).
 \hspace{3cm} $\Box$
\vspace{.5cm}

If $\bar \a = 0$, which is the equalizer case, $s_Y = Z$ and $s_X$ is undetermined.  When $\bar \a > 0$, the payoffs $s_X$ and
$s_Y$ are on the same side of $Z$, while they are on opposite sides when $\bar \a < 0$.  To be precise, we have the following.

\begin{cor}\label{newcorpressdyson}  Assume that X uses a plan $\pp$ with X Press-Dyson vector
$\tp = \g(\bar \a \SX \ + \ \bar \b \SY \ + \ \one), \ \g > 0$. Let $-Z^{-1} = \bar \a + \bar \b$. Assume that
Y uses an arbitrary general plan.

\begin{enumerate}
\item[(a)] If $\bar \a = 0$ then $s_Y = Z$.  If $\bar \a \not= 0$ then the following are equivalent
\begin{itemize}
\item[(i)] $s_Y = s_X$.
\item[(ii)] $s_Y = Z$.
\item[(iii)] $s_X = Z$.
\end{itemize}

\item[(b)] If $s_Y > s_X$ then
\begin{equation}\label{29ab}
\begin{cases}
\ \quad \bar \a > 0, \quad \ \Rightarrow \quad \ Z \ > \ s_Y \ > \ s_X. \\
\ \quad \bar \a = 0, \quad \ \Rightarrow \quad \ Z \ = \ s_Y \ > \ s_X. \\
\ \quad \bar \a < 0, \quad \ \Rightarrow \quad \ s_Y \ > \ Z \ > \ s_X.
\end{cases}
\end{equation}

\item[(c)] If $s_X > s_Y$ then
\begin{equation}\label{29aa}
\begin{cases}
\ \quad \bar \a > 0, \quad \ \Rightarrow \quad \ s_X \ > \ s_Y \ > \ Z. \\
\ \quad \bar \a = 0, \quad \ \Rightarrow \quad \ s_X \ > \ s_Y \ = \ Z. \\
\ \quad \bar \a < 0, \quad \ \Rightarrow \quad \ s_X \ > \ Z \ > \ s_Y.
\end{cases}
\end{equation}
\end{enumerate}
\end{cor}

{\bfseries Proof:} (a) If $\bar \a = 0$ then $s_Y = Z$ by (\ref{newnewpressdyson2}). If $\bar \a \not= 0$ then
(i) $\Leftrightarrow$ (ii) by (\ref{newnewpressdyson2}) and (ii) $\Leftrightarrow$ (iii) (\ref{newnewpressdyson}).

(b), (c) If $\bar \a \not= 0$ then by (\ref{newnewpressdyson2}) $s_Y - Z$ has the same sign as that of
$\bar \a (s_X - s_Y)$. \hspace{3cm} $\Box$
\vspace{.5cm}

For $Z = R$ (\ref{newnewpressdyson}) is (\ref{new16a}). When $\bar \a > 0$ these
 are the complier strategies, i.e. the generous, good plans described in Proposition \ref{prop2.4b}.

For $Z = P, \bar \a > 0$ the plans are firm.  These were considered by Press and Dyson who called
 them \emph{extortion strategies}.
The name comes from the observation that whenever Y chooses a strategy so that her  payoff is above $P$, the bonus
beyond $P$ is divided between X and Y in a ratio of $1:\kappa$. They point out that the best reply against such an extortion play by
 X is for Y is to play
$AllC = (1, 1, 1, 1)$ which gives X a payoff above $R$.   At first glance, it seems hard to escape from this coercive
effect.  I believe that the answer is for Y to play a generous good plan like the compliers above.
With repeated play, each player receives enough data to estimate statistically the strategy used by the opponent. Y's
good plan represents a credible invitation for X to switch to an agreeable plan and receive $R$, or else be locked below $R$.
Hence, it undercuts the threat from X to remain extortionate.
\vspace{.5cm}

In order to compute what happens when both players use a ZDS, we need to examine the symmetry between the
two players. Let $Switch : \R^4 \to \R^4$ be defined by $Switch(x_1,x_2,x_3,x_4) = (x_1,x_3,x_2,x_4)$. Notice that
$Switch$ interchanges the vectors $\SX$ and $\SY$. If X uses $\pp$ and Y uses
$\qq$ then recall that the response vectors  used to build the Markov matrix $\MM$ are $\pp$ and $Switch(\qq)$.
Now suppose that the two players exchange
plans so that X uses $\qq$ and Y uses $\pp$.  Then the X response is $\qq = Switch(Switch(\qq))$ and the
Y response is $Switch(\pp)$. Hence, the new Markov matrix is obtained by transposing both the second and third rows and
 the second and third columns.  It follows that if $\vv$ was a stationary vector for $\MM$, then $Switch(\vv)$
is a stationary vector for the new matrix. Hence, Theorem \ref{newtheo2} applied to the X
Press-Dyson vector $\tilde \qq$ implies that
$0 = < Switch(\vv) \cdot \tilde \qq > = < \vv \cdot Switch(\tilde \qq) >.$
Furthermore, if $\tilde \qq = a \SX + b \SY + g \one + \d \ee_{23}$,
then $Switch(\tilde \qq) = b \SX + a \SY + g \one + \d \ee_{23}$.

For a plan $\qq$, we define \emph{Y Press-Dyson vector} $\tq = Switch(\tilde \qq) = Switch(\qq) - \ee_{13}$, where
$\ee_{13} = (1, 0, 1, 0)$.  For any general plan for X and any limiting distribution $\vv$ when Y uses $\qq$ we have
$< \vv \ \cdot \ \tq > = 0$. The plan $\qq$ is a ZDS associated with $(\bar a,\bar b)$ in the ZDSstrip when
$\tq = g (\bar b \SX + \bar a \SY + \one)$ with some $g > 0$.

Now we compute what happens when X and Y use ZDS plans  associated, respectively,
with points $ (\bar \a, \bar \b)$ and
$(\bar a, \bar b)$ in the ZDS strip. This means that for some $\g > 0, g > 0$,
$\tp = \g( \bar \a S_X + \bar \b S_Y + {\mathbf 1})$ and $\tq = g( \bar b S_X + \bar a S_Y + {\mathbf 1})$.
We obtain two Press-Dyson equations which hold simultaneously
\begin{equation}\label{27}
\begin{split}
\bar \a s_X \ + \ \bar \b s_Y \quad = \quad -1, \\
\bar b s_X \ + \ \bar a s_Y \quad = \quad -1.
\end{split}
\end{equation}

If $\bar \a = \bar \b = \bar a = \bar b = -1$, which we will call  a \emph{Vertex plan} $ = \g (2(1-R), 1, 0, 1 - 2P)$,
 then the two equations are the same. Following the Remark after Lemma \ref{lem3.1},
 a Vertex plan can occur only when  $P \leq \frac{1}{2}$. Clearly, $\{ cd \}$ and $\{ dc \}$ are both terminal sets
 when both players use a Vertex plan and so the payoffs depend upon the initial plays.
  If the two players use the same initial play as well as the same plan, then $s_X = s_Y$ and the single equation
of (\ref{27}) yields  $s_X = s_Y =  \frac{1}{2}$.

Otherwise,  Lemma \ref{lem3.1} implies that
 the determinant $D \ = \ \bar \b \bar b - \bar \a \bar a$ is positive and by Cramer's Rule we get
\begin{equation}\label{28}
\begin{split}
s_X = D^{-1} (\bar a \ - \ \bar \b), \quad s_Y = D^{-1}(\bar \a - \bar b), \\
\mbox{and so} \qquad s_Y - s_X = D^{-1} [(\bar \a + \bar \b) \ - \ (\bar a + \bar b)].
\end{split}
\end{equation}
Notice that $s_X$ and $s_Y$ are independent of $ \g $ and $g$.

Thus, when both X and Y use ZDS plans from the ZDS strip, these long-term payoffs depend only
on the plans and so the results are independent of the choice of initial plays.

\begin{prop} \label{prop3.1} Assume that $\tp = \g( \bar \a S_X + \bar \b S_Y + {\mathbf 1})$ and
$\tq = g( \bar b S_X + \bar a S_Y + {\mathbf 1})$. Let $\bar \a + \bar \b = -Z_X^{-1}$ and
$\bar a + \bar b = -Z_Y^{-1}$. Assume that $(-1, -1)$ is not equal to both $(\bar \a,\bar \b)$ and $( \bar a,\bar b)$.
 \begin{itemize}

 \item[(a)] The points $(\bar \a,\bar \b), ( \bar a,\bar b)$ lie on the same value line
 $x + y = -Z^{-1}$, i.e. $Z_X = Z_Y$, iff $s_X = s_Y$.  In
that case,  $Z_X = s_X = s_Y = Z_Y$.

\item[(b)] $s_Y \ > \ s_X $ iff  $ Z_X \ > \ Z_Y $.

\item[(c)] Assume $ Z_X \ > \ Z_Y $. The following implications hold.
\begin{equation}\label{29a}
\begin{split}
\begin{cases}
\ \quad \bar \a > 0, \quad \ \Rightarrow \quad \ Z_X \ > \ s_Y \ > \ s_X. \\
\ \quad \bar \a = 0, \quad \ \Rightarrow \quad \ Z_X \ = \ s_Y \ > \ s_X. \\
\ \quad \bar \a < 0, \quad \ \Rightarrow \quad \ s_Y \ > \ Z_X \ > \ s_X.
\end{cases} \\    \\
\begin{cases}
\  \quad \bar a > 0, \quad \ \Rightarrow \quad \ s_Y \ > \ s_X \ > \ Z_Y. \\
\  \quad \bar a = 0, \quad \ \Rightarrow \quad \ s_Y \ > \ s_X \ = \ Z_Y. \\
\  \quad \bar a < 0, \quad \ \Rightarrow \quad \ s_Y \ > \ Z_Y \ > \ s_X. \\
\end{cases}
\end{split}
\end{equation}
\end{itemize}
\end{prop}

{\bfseries Proof:} We are excluding by assumption the case
when both players use Vertex plans and so we have $D > 0$.

(a) Assume $Z_X \ = \  Z_Y $. From (\ref{28}) we see that $s_Y - s_X = 0$.

When $s_X = s_Y$ Corollary \ref{newcorpressdyson}(a) implies that $s_X = s_Y = Z_X$. By using the XY symmetry
 we see that the common value is $Z_Y$ as well.
 Hence, $ Z_X = Z_Y$ and the points
lie on the same line.

(b) Since  $D > 0$, (b) follows from (\ref{28}).

(c) From (b), $ s_Y - s_X > 0 $. The first part follows from (\ref{29ab}) with $Z = Z_X$.  The second follows
from (\ref{29aa}) by using the XY symmetry with $\bar \a, \bar \b, Z$ replaced by $\bar a, \bar b, Z_Y$.
\hspace{3cm} $\Box$ \vspace{.5cm}

{\bfseries Remark:}  If both players use a Vertex plan and the same initial play, then (a) holds with
$Z_X = Z_Y = \frac{1}{2}$.

\vspace{1cm}

\section{Dynamics Among Zero Determinant Strategies}

In this section we move beyond the classical question which motivated our original interest in good strategies. We
consider now the evolutionary dynamics among memory one strategies.  We follow Hofbauer and Sigmund [9] Chapter 9 and
Akin [2].

The dynamics that we consider takes place in the context of a symmetric two-person game, but generalizing our initial
description, we merely assume that there is a  set  of strategies indexed by a finite set $\I$.
 When players X and Y use strategies
with index $i, j \in \I$,
respectively, then the payoff to player X is given by $A_{ij}$ and the payoff to Y is $A_{ji}$. Thus, the game is
described by the payoff matrix $\{ A_{ij} \}$. We imagine a population of players each using a particular strategy for each
encounter and let $\pi_i$ denote the ratio of the number of $i$ players to the total population.  The frequency
vector $\{ \pi_i \}$ lives in the unit simplex $\Delta \subset \R^{\I}$, i.e. the entries are nonnegative and sum to $1$.
The vertex $v(i)$ associated with $i \in \I$ corresponds to a population consisting entirely of $i$ players. We assume the
population is large so that we can regard $\pi$ as changing continuously in time.

Now we regard the payoff in units of \emph{fitness}.  That is, when an $i$ player meets a $j$ player in an interval of
time $dt$, the payoff $A_{ij}$ is an addition to the background reproductive rate $\rho$ of the members of the population. So
the $i$ player is replaced by $1 + (\rho + A_{ij})dt  \ i$ players. Averaging over the current population distribution,
the expected relative reproductive rate for the subpopulation of $i$ players is $\rho + A_{i \pi}$, where
\begin{equation}\label{32}
\begin{split}
A_{i \pi} \quad = \quad \Sigma_{j \in  \ \I} \ \pi_j A_{ij} \qquad \mbox{ and} \hspace{2cm}\\
A_{\pi \pi} \quad = \quad \Sigma_{i \in  \ \I} \ \pi_i A_{i \pi}  \quad = \quad \Sigma_{i,j  \ \in \ \I} \ \pi_i \pi_j A_{ij}.
\end{split}
\end{equation}

The resulting dynamical system on $\Delta$ is given by the \emph{Taylor-Jonker Game Dynamics Equations} introduced in
Taylor and Jonker [18].

\begin{equation}\label{33}
\frac{d \pi_i}{dt} \quad = \quad \pi_i (A_{i \pi} \ - \ A_{\pi \pi} ).
\end{equation}

This system is an
example of the \emph{replicator equations} studied in great detail in Hofbauer and Sigmund [9].

We will need some
general game dynamic results for later application. Fix the game matrix $\{ A_{ij} \}$.

A subset $A $ of $\Delta$ is called \emph{invariant} if $\pi(0) \in A$ implies that the entire solution
path lies in $A$.  That is, $\pi(t) \in A$ for all $t \in \R$. An invariant point is  is an
\emph{equilibrium}.

Each nonempty subset $\J$ of $ \I$ determines the \emph{face} $\Delta_{\J}$ of the simplex consisting of those $\pi \in \Delta$
such that $\pi_i = 0 $ for all $i \not\in \J$.
Each face of the simplex is  invariant  because $\pi_i = 0$ implies
that $ \frac{d \pi_i}{dt} = 0$. In particular,
 for each $i \in \I$ the vertex $v(i)$, which represents fixation at the $i$ strategy, is an equilibrium. In general,
 $\pi$ is an equilibrium when, for all $i, j \in \I$, $\pi_i, \pi_j > 0$ imply $A_{i \pi} = A_{j \pi}$.  This
 implies that $A_{i \pi} = A_{\pi \pi}$ for all $i$ such that $\pi_i > 0$.  That is, for all $i$ in the
 \emph{support} of $\pi$.

 An important example of an invariant set is the \emph{omega limit point set of an orbit}. Given an initial point
  $\pi \in \Delta$ with
 associated solution path $\pi(t)$, it is defined by intersecting the closures of the tail values.
 \begin{equation}\label{omega}
 \omega(\pi) \quad = \quad \bigcap_{t > 0} \overline{ \{ \pi(s) : s \geq t \}}.
 \end{equation}
By compactness this set is nonempty. A point is in $\omega(\pi)$ iff it is the limit of some
sequence $\{ \pi(t_n) \}$ with $\{ t_n \}$ tending to infinity.
The set  $\omega(\pi)$ consists of a single point $\pi^*$ iff $Lim_{t \to \infty} \pi(t) = \pi^*.$
In that case, $\{ \pi^* \}$ is an invariant point, i.e. an equilibrium.

\begin{df}\label{def-ess} We call a strategy $i^*$ an
\emph{evolutionarily stable strategy} (hereafter, an ESS) when
\begin{equation}\label{ess}
A_{j i^*} \ < \ A_{i^* i^*} \qquad \mbox{for all} \ \ j \not= i^* \quad \mbox{in} \ \ \I. \hspace{2cm}
\end{equation}
We call a strategy $i^*$ an
\emph{evolutionarily unstable strategy} (hereafter, an EUS) when
\begin{equation}\label{eus}
A_{j i^*} \ > \ A_{i^* i^*} \qquad \mbox{for all} \ \ j \not= i^* \quad \mbox{in} \ \ \I. \hspace{2cm}
\end{equation}
\end{df}
\vspace{.5cm}

The ESS condition above is really a special case of a more general notion, see page 63 of [9], and  is referred to there
as a \emph{strict Nash equilibrium}. We will not need the generalization and we use the term to avoid confusion with
the strategies of Nash type considered in the previous sections.

\begin{prop}\label{prop-ess} If $i^*$ is an ESS then the vertex $v(i^*)$
 is an attractor, i.e. a locally stable equilibrium, for the system (\ref{33}).  In fact, there exists $\e > 0$ such that
\begin{equation}\label{34a}
1 \ > \ \pi_{i^*} \ \geq \ 1 - \e \qquad \Longrightarrow \qquad  \frac{d \pi_{i^*}}{dt} \ > \ 0.
\end{equation}
Thus, near the equilibrium $v(i^*)$, which is characterized by $\pi_{i^*} = 1, \ \pi_{i^*}(t)$
increases monotonically, converging to $1$ and the
alternative strategies are eliminated from the population in the limit.

If $i^*$ is an EUS then the vertex $v(i^*)$
 is a repellor, i.e. a locally unstable equilibrium, for the system (\ref{33}).  In fact, there exists $\e > 0$ such that
\begin{equation}\label{34b}
1 \ > \ \pi_{i^*} \ \geq \ 1 - \e \qquad \Longrightarrow \qquad  \frac{d \pi_{i^*}}{dt} \ < \ 0.
\end{equation}
Thus, near the equilibrium $v(i^*) \ \pi_{i^*}(t)$ decreases monotonically, until the
system enters, and then remains in, the region where $\ \pi_{i^*} < 1 - \e$.
\end{prop}

{\bfseries Proof:}  When $i^*$ is an ESS, $ A_{i^* i^*} > A_{j i^*}$ for all $j \not= i^*$.
It then follows for $\e > 0$ sufficiently small that
$  \pi_{i^*}   \geq  1 - \e$ implies $A_{i^* \pi} > A_{j \pi}$ for all $j \not= i^*$. If also $1 > \pi_{i^*}, $ then
$A_{i^* \pi} > A_{\pi \pi}$. So (\ref{33}) implies (\ref{34a}).

The EUS case is similar. Notice that no solution path can cross $\Delta \cap \{ \pi_{i^*}  = 1 - \e \} $ from
$\{ \pi_{i^*}  < 1 - \e \} $.  \hspace{3cm}$\Box$
\vspace{.5cm}

\begin{df}\label{def-dom} For $\J$ a nonempty subset of $ \I $ we say a strategy $i$
\emph{weakly dominates} a strategy $j$  in $\J$ when $i, j \in \J$ and
\begin{equation}\label{dom}
A_{j k} \ \leq \ A_{i k} \qquad \mbox{for all} \ \ k  \in  \ \J, \hspace{2cm}
\end{equation}
and the inequality is strict either for $k = i$ or $k = j$.  If the inequalities are strict for all $k$ then
we say that $i$ \emph{dominates} $j$ in $\J$.

We say that $i \in \J$ dominates a sequence $\{ j_1, ..., j_n \}$ in $\J$ when $i$ dominates $j_1$ in $\J$ and
for $p = 2, ..., n$, $i$ dominates $j_p$ in $\J \setminus \{ j_1,..., j_{p - 1} \}$.

When $\J$ equals all of $\I$ we will omit the phrase ``in $\J$".
\end{df}
\vspace{.5cm}

For $i, j \in \I$, define the set $Q_{ij}$ and on it the real valued function $L_{ij}$ by
\begin{equation}\label{Lij}
\begin{split}
Q_{ij} \quad = \quad \{ \pi \in \Delta : \pi_i, \pi_j > 0 \}  \\
L_{ij}(\pi) \quad = \quad \ln(\pi_i) - \ln(\pi_j).
\end{split}
\end{equation}

\begin{lem}\label{lem-dom} (a) If $i$ weakly dominates $j$ then $dL_{ij}/dt \  > \  0$ on the
 set $Q_{ij}$.

(b) If $i$ dominates $j$ in $\J$ then there exists $\e > 0$ such that $dL_{ij}/dt \  > \  0$
on the set $Q_{ij} \cap \{ \pi \in \Delta : \Sigma_{k \not\in \J} \ \pi_k \ \leq \ \e  \}$.
\end{lem}

{\bfseries Proof:}  Observe that
\begin{equation}\label{Lijeq}
 dL_{ij}/dt \ = \ A_{i \pi} - A_{j \pi} \ = \
\Sigma_{k \in  \ \I} \pi_k (A_{ik} - A_{jk})
\end{equation}

(a) Since $\pi_i, \pi_j > 0$ in $Q_{ij}$
and $A_{ik} - A_{jk} \geq 0$ for all $k$ with strict inequality for
$k = i$ or $k = j$, it follows that the derivative is positive. \hspace{3cm}

(b)  Define
\begin{equation}\label{manydef}
\begin{split}
m \quad = \quad min \{ A_{ik} - A_{jk} : k \in \J \} \ > \ 0, \\
M \quad = \quad  max \{ |A_{ik} - A_{jk}| : k \not\in \J \},\\
\pi_{\J} \quad = \quad \Sigma_{k \in  \ \J} \pi_k, \hspace{2cm} \\
\pi_{k|\J} \quad = \quad \pi_k / \pi_{\J} \quad \mbox{for} \ k \in \J. \hspace{1cm}
\end{split}
\end{equation}

Observe that $\Sigma_{k \not\in  \ \J} \pi_k  = 1 - \pi_{\J}$.

For any $\pi \in Q_{ij}$
\begin{equation}
\begin{split}
A_{i\pi} - A_{j\pi} \ = \   \pi_{\J}  \Sigma_{k \in  \ \J} \ \pi_{k|\J} (A_{ik} - A_{jk}) \
+ \ \Sigma_{k \not\in  \ \J} \pi_k (A_{ik} - A_{jk}) \\
\geq \quad  \pi_{\J} m   \ - \ (1 - \pi_{\J}) M. \hspace{3cm}
\end{split}
\end{equation}
So if $\e$ is chosen with $0 < \e < m/(m + M) $ then $A_{i\pi} - A_{j\pi} > 0$
when $\pi \in Q_{ij} \cap \{ \pi \in \Delta : (1 - \pi_{\J}) \  \leq \ \e  \}$. \hspace{3cm}
$\Box$ \vspace{.5cm}

\begin{lem}\label{lem-lim}
 If $\pi(t)$ is a solution path with $\pi(0) \in Q_{ij}$ and there exists $T \in \R$ such
 that $dL_{ij}/dt \  > \  0$
on the set $Q_{ij} \cap \overline{\{ \pi(t) : t \geq T \}}$, then
\begin{equation}\label{j-dom}
Lim_{t \to \infty} \ \pi_j(t) \quad = \quad 0.  \hspace{2cm}
\end{equation}
\end{lem}

{\bfseries Proof:}  By assumption, $ L_{ij}(\pi(t))$ is a strictly increasing function of $t$
for $t \geq T$. Thus, as a $t$ tends to infinity $L_{ij}(\pi(t))$ approaches
$\ell = sup \{ L_{ij}(\pi(t)) : t \geq T \}$ with
$L_{ij}(\pi(T)) < \ell \leq + \infty$.

We must prove  that $\pi_j = 0$ on the omega limit set.
Assume instead that $\pi^* \in \omega(\pi(0))$ with $\pi^*_j > 0$.  If $\pi^*_i$ were $0$ then $L_{ij}(\pi(t))$ would
not be bounded below on $\{ \pi(t) : t \geq T \}$.  Hence, $\pi^*$  lies in
$Q_{ij}$ with $\ell = L_{ij}(\pi^*) < \infty$. So on the invariant set
$\omega(\pi(0)) \cap Q_{ij}$, which contains $\pi^*$ and so is nonempty, $L_{ij}$ would be constantly $\ell < \infty$. Since
this set is invariant, $ dL_{ij}/dt $ would equal zero. This contradicts
our assumption that the derivative is positive on $\omega(\pi(0)) \cap Q_{ij}$. \hspace{3cm} $\Box$
\vspace{.5cm}

\begin{prop}\label{prop-dom} For $i \in \I$, let $\pi(t)$ be a solution path with $\pi_i(0)  > 0$

(a) If $i$ weakly dominates $j$ then $Lim_{t \to \infty} \ \pi_j(t) \  = \  0.$

(b) If $i$ dominates the sequence $\{ j_1,...,j_n \}$ then for
$j = j_1, ..., j_n, \ Lim_{t \to \infty} \ \pi_j(t) \quad = \quad 0. $
\end{prop}

{\bfseries Proof:} (a) If $\pi_j(0) = 0$, then $\pi_j(t) = 0$ for all $t$ and so the limit is $0$. Hence,
we may assume $\pi_j(0) > 0$ and so that $\pi(0) \in Q_{ij}$.  By Lemma \ref{lem-dom} (a), $dL_{ij}/dt \  > \  0$ on
 $Q_{ij} $ and so Lemma \ref{lem-lim} implies $Lim_{t \to \infty} \ \pi_j(t) \  = \  0. $

 (b) We prove the result by induction on $n$.

 By part (a) $Lim_{t \to \infty} \ \pi_j(t) \  = \  0. $  for $j = j_1$.

 Now assume the limit result is true
 for $j =  j_1,...,j_{p-1} $ with $1 < p \leq n$.  We prove the result for $j = j_p$.

 Let $\J = \I \setminus \{ j_1,...,j_{p-1}  \}$ . By assumption, $i $ dominates $j_p$ in $\J$. Hence, with $j = j_p$
 Lemma \ref{lem-dom} (b) implies there exists $\e > 0$ such that    $dL_{ij}/dt \  > \  0$
on the set $Q_{ij} \cap \{ \pi \in \Delta : \Sigma_{k \not\in \J} \pi_k \leq \e  \}$.

By induction hypothesis, there exists $T$ such that $\Sigma_{k \not\in \ \J} \pi_k(t) \leq \e $ for all $t \geq T$.
Hence, $\overline{\{ \pi(t) : t \geq T \} } \subset \{ \pi : \Sigma_{k \not\in \J} \pi_k(t) \leq \e\}$.

As in part (a), we can assume $\pi \in Q_{ij}$ and then apply
Lemma \ref{lem-lim} to conclude $Lim_{t \to \infty} \ \pi_j(t) \  = \  0. $  This completes the inductive step.
 \hspace{1cm} $\Box$
 \vspace{.5cm}

Now we specialize to the Iterated Prisoner's Dilemma.
By a \emph{strategy} we will mean a  plan $\pp$ together with an initial play, pure or mixed. Recall that a
good (or agreeable) strategy is a good (resp. agreeable) plan together with initial cooperation.

To apply the Taylor-Jonker dynamics to our case, we suppose that $\I$ indexes  a finite collection
of strategies.
We then use
\begin{equation}\label{34}
A_{ij} \quad = \quad s_X \qquad \mbox{so that} \qquad A_{ji} \quad = \quad s_Y.
\end{equation}
That is, when  the X player uses the $i$ strategy and the Y player uses the $j$ strategy then the players
receive the payoffs $s_X$ and $s_Y$, respectively, as additions to their reproductive rate.
When the associated  Markov matrix is convergent,
there is a unique terminal set, and  the long term payoffs, $s_X, s_Y$ depend only on the plans
and not on the initial plays.

\begin{theo}\label{thm-ess-good} Let $\I$ index a finite set of strategies for the Iterated
Prisoner's Dilemma. Suppose that associated with $i^* \in \I$
is a good strategy $\pp^{i*}$. If for no other $j \in \I$ is the plan
$\pp^j$ agreeable, then $i^*$ is an ESS for the associated game $\{ A_{ij} : i,j \in \I \}$ and so the vertex $v(i^*)$
is an attractor for the dynamic. \end{theo}

{\bfseries Proof:} Since $i^*$ is associated with an agreeable strategy, $A_{i^* i^*} = R$.
Since $\pp^{i^*}$ is good and  $\pp^j$ is not agreeable for $j \not= i^*$, it follows from Corollary \ref{Hilbecor}(a) that
$A_{j i^*} < R$ for $j \not= i^*$. Thus, $i^*$ is an ESS.  $\Box$
\vspace{.5cm}

There are other cases of ESS which are far from good.

\begin{lem}\label{lem-ess-bad} (a) Assume that X uses a plan $\pp = (p_1, p_2, 0, 0)$ with $p_1, p_2 < 1$.  If Y
uses any plan $\qq$ which is not firm, then
\begin{equation}\label{ess-bad}
s_Y \ < \ P \ < s_X. \hspace{3cm}
\end{equation}

(b) Assume that X uses a plan $\pp = (1, p_2, 0, 0)$ with $ p_2 < 1$.  If Y
uses any plan $\qq$ which is neither firm nor agreeable, then  (\ref{ess-bad}) holds.

(c) Assume $P < \frac{1}{2}$ and that X uses a firm, non-exceptional ZDS with $\bar \a < 0$.   If Y
uses any plan $\qq$ which is not firm, then  (\ref{ess-bad}) holds.
\end{lem}

{\bfseries Proof:} (a) and (b) Since $p_3 = p_4 = 0$ the set $\{ dc, dd \}$ is closed. If $q_4 > 0$ then $\{ dd \}$ is not
closed and so is not a terminal set.

(a) Since $p_2 < 1$ there is an edge from $cd$ to either $dc$ or $dd$.  Hence, $cd$ is transient. Similarly,
$p_1 < 1$ implies $cc$ is transient. Hence, for any stationary distribution $\vv$, $v_1 = v_2 = 0$.
Since $\qq$ is not firm, $q_4 > 0$ and so $v_4 < 1$.   Hence, $s_Y = v_4 P < P$ and
$s_X = v_3 + v_4 P  = (1 - v_4) + v_4 P > P$.

(b) As before $p_2 < 1$ implies that $cd$ is transient. Now $\qq$ is not agreeable and so $q_1 < 1$. This implies
there is an edge from $cc$ to the transient state $cd$ and so $cc$ is transient. The proof is completed as in (a).

(c)  Because $P < \frac{1}{2}$, the smallest entry in $\frac{1}{2}(\SX + \SY)$ is P and so
$\frac{1}{2}(s_X + s_Y) \leq P$  can only happen when $v_4 = 1$ which implies $s_X = s_Y = P$. This  requires that Y play a
firm plan so that $\{ dd \}$ is a terminal set.  Compare Proposition \ref{newprop1}.

From  (\ref{newnewpressdyson}) we see that with $\bar \a + \bar \b = - Z^{-1}$ and $\kappa = \bar \a Z/(1 + \bar \a Z)$
\begin{equation}\label{newnewpressdyson3}
\frac{1}{2}(1 + \kappa)(s_X \ - \ Z) \quad = \quad (\frac{1}{2}(s_X + s_Y) \ - \ Z).
\end{equation}
When $Z = P$,
 $P < \frac{1}{2}$ and $-1 \geq \bar \a$ imply that $(1 + \kappa) = (1 + 2 \bar \a P)/(1 + \bar \a P) > 0$. Hence,
 $s_X \leq P$ implies $\frac{1}{2}(s_X + s_Y) \leq P$.  Since the Y plan is not firm, this does not happen. Hence,
 $s_X > P$. Since $\kappa < 0$, (\ref{newnewpressdyson}) implies that $s_Y < P$.
\hspace{2cm} $\Box$
\vspace{.5cm}

\begin{theo}\label{thm-ess-bad} Let $\I$ index a finite set of strategies for the Iterated Prisoners Dilemma.

(a) Suppose that associated with $i^* \in \I$
is a  plan $\pp^{i*} = (p_1, p_2, 0, 0)$ with $p_1, p_2 < 1$ together with any initial play.
If for no other $j \in \I$ is the plan
$\pp^j$ firm, then $i^*$ is an ESS for the associated game $\{ A_{ij} : i,j \in \I \}$.

(b) Suppose that associated with $i^* \in \I$
is a  plan $\pp^{i*} = (1, p_2, 0, 0)$ with $ p_2 < 1$ together with any initial play.
If for no other $j \in \I$ is the plan
$\pp^j$ either agreeable or firm, then $i^*$ is an ESS for the associated game $\{ A_{ij} : i,j \in \I \}$.

(c) Assume that $P < \frac{1}{2}$.
Suppose that associated with $i^* \in \I$ is a firm, non-exceptional ZDS with $\bar \a < 0$ together with
any initial play. If for no other $j \in \I$ is the plan
$\pp^j$ firm, then $i^*$ is an ESS for the associated game $\{ A_{ij} : i,j \in \I \}$.
\end{theo}

{\bfseries Proof:} (a) If both players use $p^{i^*}$ then there is an edge from $dc$ to $dd$ and so
$dc$, $cd$ and $cc$ are all transient.  Thus, $\{ dd \}$ is the unique terminal set and so $A_{i^* i^*} = P$ regardless
of the initial plays. By Lemma \ref{lem-ess-bad}(a)
$A_{j i^*} < P$ for all $j \not= i^*$.

(b) If both players use $p^{i^*}$ then there are edges from $cd$ to $dd$ and from $dc$ to $dd$. The two
terminal sets are $\{ cc \}$ and $\{ dd \}$.  Hence, $ R \geq A_{i^* i^*} \geq P$.
This time Lemma \ref{lem-ess-bad}(b) implies
$ A_{j i^*} < P $ for any $j \not= i^*$.

(c) $A_{i^* i^*} = Z_{i^*} = P$ since the $i^*$ plan is firm. Lemma \ref{lem-ess-bad}(c) implies
$ A_{j i^*} < P $ for any $j \not= i^*$.
\hspace{3cm} $\Box$
\vspace{.5cm}

Thus, $p^{i^*} = All D = (0, 0, 0, 0)$ with any initial play
is an ESS when played against plans which are not firm.
If $p^{i^*} = Grim = (1, 0, 0, 0)$ then with any initial play $i^*$ is an ESS
when played against strategies which are neither agreeable nor firm.

At the other extreme we have the following.

\begin{theo}\label{thm-eus} Let $\I$ index a finite set of strategies for the Iterated Prisoners' Dilemma.
  Assume that $P < \frac{1}{2}$.
 Suppose that associated with $i^* \in \I$
is an extortionate plan $\pp^{i*}$ together with initial defection.  That is, $\pp^{i^*}$ is a firm
ZDS with $\bar \a > 0$. If for no other $j \in \I$ is the plan
$\pp^j$ firm, then $i^*$ is an EUS for the associated game $\{ A_{ij} : i,j \in \I \}$ and so the vertex $v(i^*)$
is a repellor for the dynamic. \end{theo}

{\bfseries Proof:} Because $P < \frac{1}{2}$, the smallest entry in $\frac{1}{2}(\SX + \SY)$ is P and so
$s_X,s_Y \leq P$ implies $s_X = s_Y = P$ and this can only happen when $v_4 = 1$ which requires that Y play a
firm plan so that $\{ dd \}$ is a terminal set.  Compare Proposition \ref{newprop1}.

Since the $i^*$ strategy is firm with initial defection, $A_{i^* i^*} = P$.

If Y uses any plan which is not firm then (\ref{newnewpressdyson}) with $z = P$ and $\bar \a > 0$ shows that
if $s_Y \leq P$ then $s_X \leq P$ as well.  Because $P < \frac{1}{2}$ this can only happen when $s_X = s_Y = P$ and
$v_4 = 1$.  But the Y plan is not firm.  It follows that $s_Y > P$. Thus, for any $j \not= i^*$, $A_{ji^*} > A_{i^* i^*}$.
This says that strategy $i^*$ is an EUS. \hspace{3cm}$\Box$
\vspace{.5cm}

We now specialize to the case when all the strategies indexed by $\I$ are ZDS's with
the exceptional strategies excluded. We can thus regard
$\I$ as listing a finite set of points $(\bar \a_i,\bar \b_i)$ in the ZDSstrip and, except for the Vertex plans,
we may disregard the initial plays.
We define $Z_i = - (\bar \a_i + \bar \b_i)^{-1}$. That is, the point $(\bar \a_i,\bar \b_i)$
lies on the value line $x + y = - (Z_i)^{-1}$.

 X uses $\pp$ associated with
$(\bar \a_i, \bar \b_i)$ when $\tp = \g_i( \bar \a_i \SX + \bar \b_i \SY + {\mathbf 1})$ and Y uses $\qq$
associated with $(\bar \a_j, \bar \b_j)$ when $\tq = \g_j( \bar \b_j \SX + \bar \a_j \SY + {\mathbf 1})$ for some
$\g_i, \g_j > 0$. Notice the XY switch.

If both players use a Vertex plans with the same initial plays, then
$(\bar \a_i, \bar \b_i)  = (-1,-1)$ and $A_{ii} = \frac{1}{2} = Z_i$.
Recall that $(-1,-1)$ lies in the ZDSstrip iff $P \leq \frac{1}{2}$.

 Otherwise, we apply (\ref{27})
with $(\bar \a,\bar \b) = (\bar \a_i,\bar \b_i)$ and $(\bar a,\bar b) = (\bar \a_j,\bar \b_j)$. Then from
(\ref{28}) we get, for $i \not= j$
\begin{equation}\label{35}
\begin{split}
A_{ij} \quad = \quad s_X \quad = \quad  K_{ij} (\bar \a_j \ - \ \bar \b_i)  \hspace{2cm}        \\
\mbox{ with} \qquad K_{ij} \ = \ K_{ji} \ = \ (\bar \b_i \bar \b_j - \bar \a_i \bar \a_j)^{-1} \ > \ 0.
\end{split}
\end{equation}
Note that these payoffs are independent of the choice of $\g_i, \g_j$ as well as the initial plays.

By  Proposition \ref{prop3.1}(a)
\begin{equation}\label{35aa}
A_{ii} \quad = \quad Z_i \qquad \mbox{for all} \ i \in \I.\hspace{3cm}
\end{equation}

We begin with some degenerate cases. For convenience, we exclude the Vertex plans.

First, if all of the points $(\bar \a_i,\bar \b_i)$ lie on the same value line $x + y = -Z^{-1}$, i.e. all the $Z_i$'s are equal,
then by Proposition \ref{prop3.1} (a)
$A_{ij} = Z$ for all $i,j$ and so $\frac{d \pi}{dt} = 0$ and every population distribution is an equilibrium.  In general, if
for two strategies $i,j \ \ A_{ij} = A_{ji} = Z$ then
by Proposition \ref{prop3.1}(a) both points lie on $x + y = -Z^{-1}$ and it follows that
$A_{ii} = A_{jj} = Z$ as well. In general, if $\I_Z = \{ i : Z_i = Z \}$ contains more than one $i \in I$ then the dynamics
is degenerate on the face $\Delta_{\I_Z}$ of the simplex.

Second, if all of the points satisfy $\bar \a_i = 0$ then all the strategies are equalizer strategies. In this case
the payoff matrix need not be constant but $A_{ij}$ depends only on $j$.  This implies that for all $i$
$A_{i \pi} = A_{\pi \pi}$  and so again $\frac{d \pi}{dt} = 0$ and every population distribution is an equilibrium.

We will now see that the line $\bar \a = 0$ separates different interesting dynamic behaviors.

\begin{theo}\label{theo3.2} Let $\I$ index a set of non-exceptional ZDS plans. Thus,   each $i \in \I$ is associated
with a point $(\bar \a_i, \bar \b_i)$ in the ZDS strip and $\bar \a_i + \bar \b_i = -( Z_i)^{-1}$.
\vspace{.25cm}

 Assume  either

Case (+): $\bar \a_i > 0$ for all $i \in \I$ and for some $i^* \in \I$, $Z_{i^*} > Z_j$ for all $j \not= i^*$;

or

Case (-): $\bar \a_i < 0$ for all $i \in \I$ and for some $i^* \in \I$, $Z_{i^*} < Z_j$ for all $j \not= i^*$.
\vspace{.25cm}

The strategy $i^*$ is an ESS and if $\pi_{i^*}(0) > 0$ then the solution path converges to the vertex $v(i^*)$.
\end{theo}

{\bfseries Proof:} List the strategies $j_1,....,j_n$ of $\I \setminus \{ i^* \}$ so that
in Case(+) $Z_{j_1} \leq Z_{j_2} \leq ...\leq Z_{j_n} < Z_{i^*}$ and in
Case(-) $Z_{j_1} \geq Z_{j_2} \geq ...\geq Z_{j_n} > Z_{i^*}$.
For both cases we apply Proposition \ref{prop3.1}. It first implies that if $Z_i = Z_j$ then
\begin{equation}\label{+and-}
A_{ii} \ = \ Z_i \ = \ A_{ji} \ = \ A_{ij} \ = \ Z_j \ = \ A_{jj}.
\end{equation}

Case(+) If $Z_i > Z_j$, then, because $ \bar \a_i, \bar \a_j > 0$, Proposition \ref{prop3.1} implies that
\begin{equation}\label{case+}
A_{ii} \ = \ Z_i \ > \ A_{ji} \ > \ A_{ij} \ > \ Z_j \ = \ A_{jj}.
\end{equation}

Hence, if $Z_i > Z_k \geq Z_j$ then $A_{ii} > A_{ji}$ and $A_{ik} > A_{kk} \geq A_{jk}$.

It follows that $i^*$ dominates the sequence $\{ j_1,...,j_n \}$. Hence, Proposition \ref{prop-dom} (b) implies
that $Lim_{t \to \infty} \pi_j(t) = 0$ for $j = j_1,...,j_n$ when $\pi_{i^*}(0) > 0$. Consequently,
$\pi_{i^*}(t) = 1 - \Sigma_{p=1}^n \ \pi_{j_p}(t)$
tends to $1$.  That is, $\pi(t)$ converges to $v(i^*)$.

Case(-)  If $Z_i < Z_j$, then, because $ \bar \a_i, \bar \a_j < 0$, Proposition \ref{prop3.1} implies that
\begin{equation}\label{case-}
A_{ij} \ > \  \ Z_j \ = \ A_{jj}  >  \ A_{ii} \ =  \ Z_i \  > \  A_{ji}.
\end{equation}

It again follows that $i^*$ dominates the sequence $\{ j_1,...,j_n \}$ and convergence to $v(i^*)$ again follows from
Proposition \ref{prop3.1}.

In both cases, it is clear that $i^*$ is an ESS. \hspace{3cm} $\Box$
\vspace{.5cm}

Thus, when only $\bar \a > 0$ ZDS plans are competing with one another, the ones on the highest value line
win.  Among $\bar \a < 0$ ZDS plans the ones on the lowest value line win.

The local stability of an ESS good strategy  will  not be global when both signs occur.  To illustrate this, consider the
case of two strategies indexed by $\I = \{ 1,2 \}$. Letting $w = \pi_1$, it is an easy exercise to show that (\ref{33}) reduces to
\begin{equation}\label{40a}
\frac{d w}{d t} \quad = \quad w (1 - w)[(A_{11} - A_{21}) w + (A_{12} - A_{22})(1 - w)].
\end{equation}

\begin{prop}\label{prop32a} Assume that $Z_1 \ > \ Z_2$ and that
$\bar \a_1 \cdot \bar \a_2 \ < \ 0$. There is an equilibrium population $\pi^* = (w^*,(1-w^*))$ which contains both strategies.
with
\begin{equation}\label{41a}
w^*/(1 - w^*) \quad = \quad (A_{22} - A_{12})/(A_{11} - A_{21}).
\end{equation}
This equilibrium is stable if $\bar \a_1 < 0$ and is unstable if $\bar \a_1 > 0$.
\end{prop}

{\bfseries Proof:}  If $\bar \a_1 < 0$ and $\bar \a_2 > 0$ then Proposition \ref{prop3.1} implies
that   $ A_{11} - A_{21} = Z_1 - A_{21} < 0 $ and   $A_{12} - A_{22} = A_{12} - Z_2 > 0$. Reversing the signs reverses the
inequalities. The result then easily follows from equation (\ref{40a}). Just graph the linear
function of $w$ in the brackets and observe where the result is positive or negative.
$\Box$ \vspace{.5cm}

\begin{ques}\label{ques} Suppose we restrict to the case where $\I$ indexes ZDS's lying on different value lines
 to avoid degeneracies. We ask: \end{ques}
\begin{itemize}
\item How large a population can coexist?  If $N$ is the size of $\I$, the number of competing strategies, then
for what  $N$ do there exist examples with an interior
equilibrium, that is, an equilibrium $\pi$ such that $\pi_i > 0$ for all $i \in \I$?  When is there a locally
stable interior equilibrium?  For how large an $N$ can \emph{permanence} occur (see [9] Section 3),
that is, the boundary of $\Delta$ be a repellor? The Brouwer Fixed Point Theorem implies that such a permanent
system always admits an interior equilibrium.  When an interior equilibrium does not exist there is always some sort
of dominance among the mixed strategies of the game $\{ A_{ij} \}$.  See [1] and [3].

\item Can there exist a stable, closed invariant set containing no equilibria, e.g. a stable limit cycle?
\end{itemize}
\vspace{.5cm}

There is alternative version of the dynamics which explicitly considers for X not the payoff $s_X$ but
the advantage that X has over Y. That is, the addition to the growth rate is given not by
$s_X$ but by the difference $s_X - s_Y$.  This amounts to replacing $A_{ij}$ by the anti-symmetric matrix
$S_{ij} = A_{ij} - A_{ji}$ so that the game becomes zero-sum. In this case, we define $\xi_i =  -Z_i^{-1} = \bar \a_i + \bar \b_i$.
Thus, $\xi_i $ varies in the interval $[-P^{-1}, -R^{-1}]$.
Define $\xi_{\pi} = \Sigma_{i  \in \ \I} \ \pi_i \xi_i$. From (\ref{35}) we get
\begin{equation}\label{39}
S_{ij} \quad = \quad K_{ij} (\xi_j \ - \ \xi_i),
\end{equation}
where we let $K_{ii} = 1$ for all $i \in \I$.

Since $\{ S_{ij} \}$ is antisymmetric, $S_{\pi \pi} \ = \ 0 $.

For this system the behavior is always like the $\bar \a < 0$ case for the previous system.

\begin{theo}\label{theo3.3} Let $\I$ index a finite list of non-exceptional ZDS strategies, with at most one using a
Vertex plan.
 For the system with
\begin{equation}\label{40}
\frac{d \pi_i}{dt} \quad = \quad \pi_i ( S_{i \pi} - S_{\pi \pi}) \quad = \quad \pi_i  S_{i \pi},
\end{equation}
we have
\begin{equation}\label{41}
\begin{split}
\frac{d \xi_{\pi}}{dt} \quad \leq \quad 0, \hspace{4cm}\\
 \mbox{with equality iff} \qquad \pi_i, \pi_j > 0 \quad \Longrightarrow \quad \xi_i \ = \ \xi_j.
 \end{split}
 \end{equation}

 Assume now that  $Z_{i^*} < Z_j$, or equivalently $\xi_{i^*} < \xi_j$  for all $j \not= i^*$. The strategy $i^*$ is an ESS and
 if $\pi_{i^*}(0) > 0$ then the solution path converges to the vertex $v(i^*)$.
 \end{theo}

 {\bfseries Proof:} Because $K_{ij}$ is symmetric and positive, $\frac{d \xi_{\pi}}{dt}$ equals
 \begin{equation}\label{42}
 \begin{split}
- \Sigma_{i,j \ \in \ \I} \ \pi_i \pi_j K_{ij} \xi_i (\xi_i - \xi_j) \quad = \hspace{2cm}\\
- \frac{1}{2}[  \Sigma_{i,j \ \in \ \I} \ \pi_i \pi_j K_{ij} \xi_i (\xi_i - \xi_j) -
 \Sigma_{i,j \ \in \ \I} \ \pi_j \pi_j K_{ij} \xi_j (\xi_i - \xi_j) ] \\
 = \quad - \frac{1}{2} \Sigma_{i,j \ \in \ \I} \  \pi_i \pi_j K_{ij} (\xi_i - \xi_j)^2 \quad \leq \quad 0.\hspace{2cm}
 \end{split}
 \end{equation}
Equality holds iff $\pi_i \pi_j  (\xi_i - \xi_j)^2 = 0$ for all $i, j \in \I$. That is, when $\xi_i = \xi_j$
for all $i, j$ with $\pi_i, \pi_j > 0$.

If $\xi_i < \xi_j$ then

\begin{equation}\label{42a}
S_{ij} \ > \  \ 0 \ = \ S_{jj}  =  \ S_{ii} \  > \  S_{ji}.
\end{equation}

If $\xi_i < \xi_k \leq \xi_j$,  then $S_{ik} > S_{kk} \geq S_{jk}$.  Let $\{ j_1,...,j_n \}$  list $\I \setminus \{ i^* \}$
with $\xi_{j_1} \geq .... \geq \xi_{j_n}$. As in Case(-) of Theorem \ref{theo3.2} it
follows that $i^*$ dominates the sequence $\{ j_1,...,j_n \}$.
If $\pi_{i^*}(0) > 0$ then $\pi(t)$ converges to $v(i^*)$ by
Proposition \ref{prop3.1}. \hspace{3cm}  $\Box$
\vspace{1cm}

\section*{References}

\begin{itemize}

\item[1.] E. Akin,  (1980) Domination or equilibrium, Math. Biosciences {\bfseries 50:} 239-250.

\item[2.] ---, (1990) The differential geometry of population genetics and evolutionary games, In S. Lessard (ed.)
{\bfseries Mathematical and statistical developments of evolutionary theory},  Kluwer, Dordrecht : 1-93.

\item[3.] E. Akin, J. Hofbauer, (1982) \emph{Recurrence of the unfit}, Math. Biosciences {\bfseries 61:} 51-63.

\item[4.] R. Axelrod, \emph{ The Evolution of Cooperation}, Basic Books, New York, NY, 1984.

\item[5.] M. Boerlijst, M. Nowak,  K. Sigmund,  (1997) Equal pay for all prisoners, Amer. Math. Monthly, {\bfseries  104}  303-305.

\item[6.] M. D. Davis, \emph{ Game Theory: A Nontechnical Introduction}, Dover Publications, Mineola, NY, 1983.

\item[7.] G. Hardin, (1968)  The tragedy of the commons, Science, {\bfseries  162}  1243-1248.

\item[8.] C. Hilbe, M. Nowak,  K. Sigmund, (2013)
The evolution of extortion in iterated Prisoner's Dilemma games,  PNAS, {\bfseries 110} no. 17,  6913-6918.

\item[9.] J. Hofbauer,  K. Sigmund,  \emph{ Evolutionary Games and Population Dynamics},
Cambridge Univ. Press, Cambridge, UK, 1998.

\item[10.] S. Karlin,  H. Taylor,  \emph{ A First Course in Stochastic Processes}, second edition, Academic Press, New York, NY, 1975.

\item[11.] J. Maynard Smith, \emph{ Evolution and the Theory of Games}, Cambridge Univ. Press, Cambridge, UK, 1982.

\item[12.] M. Nowak,  \emph{ Evolutionary Dynamics}, Harvard Univ. Press, Cambridge, MA, 2006.

\item[13.] K. Sigmund,  \emph{ Games of life}, Oxford Univ. Press, Oxford, UK, 1993.

\item[14.] ---,  \emph{\ The Calculus of Selfishness}, Princeton Univ. Press, Princeton, NJ, 2010.

\item[15.] W. Press,   F. Dyson,  (2012) Iterated Prisoner's Dilemma contains
strategies that dominate any evolutionary opponent, PNAS,   {\bfseries 109} no. 26, 10409-10413.

\item[16.] A. Stewart, J. Plotkin,  (2012) Extortion and cooperation in the Prisoner's Dilemma, PNAS, {\bfseries 109} no. 26, 10134-10135.

\item[17.] P. D. Straffin, \emph{ Game Theory and Strategy}, Mathematical Association of America, Washington, DC, 1993.

\item[18.] P. Taylor, L. Jonker,  (1978) Evolutionarily stable strategies and game dynamics, Math. Biosciences {\bfseries 40:}
145-156.

\item[19.] J. Von Neumann, O. Morgenstern,   \emph{ Theory of Games and Economic Behavior}, Princeton Univ. Press,
Princeton, NJ, 1944.

\end{itemize}

\end{document}